\documentclass[%
jmp,
 amsmath,amssymb,amsfonts,
Conference Proceedings
]{revtex4-1}

\usepackage{graphicx, amsfonts, amssymb}
\usepackage{dcolumn}
\usepackage{bm}

\usepackage[utf8]{inputenc}
\usepackage[T1]{fontenc}
\usepackage{mathptmx}
\usepackage{etoolbox}

\makeatletter
\def\@email#1#2{%
 \endgroup
 \patchcmd{\titleblock@produce}
  {\frontmatter@RRAPformat}
  {\frontmatter@RRAPformat{\produce@RRAP{*#1\href{mailto:#2}{#2}}}\frontmatter@RRAPformat}
  {}{}
}%
\makeatother

\usepackage{xcolor,ulem,textcomp,MnSymbol}

\newtheorem{thm}{Theorem}[section]
\newtheorem{Prop}[thm]{Proposition}
\newtheorem{Lem}[thm]{Lemma}
\newtheorem{Cor}[thm]{Corollary}
\newtheorem{Def}[thm]{Definition}

\newtheorem{Rmk}[thm]{Remark}

\def\TT{{\large \Xi}}

\def\R {\mathbb{R}}

\def\bbC {\mathbb{C}}
\def\bbD {\mathbb {D}}
\def\bbE {\mathbb {E}}

\def\bbN {\mathbb{N}}
\def\bbP {\mathbb {P}}
\def\bbR {\mathbb{R}}
\def\bbS {\mathbb{S}}
\def\bbT{{\mathbb T}}

\def\cA {\mathcal{A}}
\def\cB {\mathcal{B}}
\def\cC {\mathcal{C}}
\def\cE {\mathcal{E}}
 
\def\cD {\mathcal{D}}

\def\cE {\mathcal{E}}

\def\cG {\mathcal{G}}

\def\cM {\mathcal{M}}
\def\cN {\mathcal{N}}

\def\cP {\mathcal{P}}

\def\cZ {\mathcal{Z}}

\def\cR {\mathcal{R}}
\def\cT {\mathcal{T}}

\def\gL{\Lambda}

\def\gp {{\varphi}}

\def\eps {{\varepsilon}}
\def\e {{\varepsilon}}

\def\gd {{\delta}}
\def\eps {{\varepsilon}}
\def\e {{\varepsilon}}
\def\gl{\lambda}

\def\indc {{\textbf 1}}

\def\mes {\pmb{\nu}_{\eps,T}^{(H)}}
\def\mesflat {\pmb{\nu}_{\eps,T}^{(0)}}

\def\mesh0 {\pmb{\nu}_{\eps,T}^{(H),\hbox{\rm\scriptsize min}}}

\def\tanakae{  \pmb{\Upsilon}_{\eps,T}  }
\newcommand{\tanaka}[1]{\pmb{\Upsilon}_{#1}}

\begin{document}


\title[Cluster expansion for a dilute hard sphere gas dynamics]{Cluster expansion for a dilute hard sphere gas dynamics}

\author{Thierry Bodineau}
\affiliation{CMAP, CNRS, Ecole Polytechnique, I.P. Paris
\\
Route de Saclay, 91128 Palaiseau Cedex, FRANCE}
\email{thierry.bodineau@polytechnique.edu}

\author{Isabelle Gallagher}
\affiliation{DMA, \'Ecole normale sup\'erieure, CNRS, PSL Research University\\
45 rue d'Ulm, 75005 Paris, FRANCE \\
and Universit\'e de Paris}
\email{gallagher@math.ens.fr}

\author{Laure Saint-Raymond}
\affiliation{IHES, 35 Rte de Chartres, 91440 Bures-sur-Yvette,
FRANCE
}
\email{laure@ihes.fr}

\author{Sergio Simonella}
\affiliation{UMPA UMR 5669 du CNRS, ENS de Lyon,\\
46  all\'ee d'Italie, 69007 Lyon,
FRANCE \\
}
\email{sergio.simonella@ens-lyon.fr}

%

\date{\today}

\begin{abstract}
In  \cite{BGSS2}, a cluster expansion method has been developed  to study 
the fluctuations of the hard sphere dynamics  around the Boltzmann equation.
This method provides a precise control on the exponential moments of the empirical measure,
from which the fluctuating Boltzmann equation and large deviation estimates have been deduced.
The cluster expansion in  \cite{BGSS2} was implemented at the level of the BBGKY hierarchy, 
which is a standard tool to investigate the deterministic dynamics \cite{Cercignani_Illner_Pulvirenti}.
In this paper, we introduce an alternative approach, in which the cluster expansion is applied directly on 
real trajectories of the particle system. 
This offers a fresh perspective on the study of the hard sphere dynamics in the low density limit, allowing
to recover the results obtained in  \cite{BGSS2}, and also to describe the actual clustering of particle trajectories.
  
\end{abstract}

\maketitle

\section{Introduction and presentation of the model}

A gas  can be modelled by a billiard made of hard spheres, moving in agreement with the laws of classical mechanics. Initially (at time $t=0$) the spheres are randomly and identically distributed according to a probability measure which is then transported by the flow of the deterministic dynamics (see~(\ref{hardspheres})-(\ref{defZ'nij}) below). For a dilute gas, it has been shown in the seminal 
work of Lanford \cite{Lanford} that when the average number of particles goes to infinity in the Boltzmann-Grad limit, the gas density converges towards a solution to the Boltzmann equation, at least for a short time. 
This work triggered a wave of developments, including some recent quantitative convergence results and generalisations to the case of compactly supported potentials; see \cite{spohn2012largeO, Cercignani_Illner_Pulvirenti, Cercignani_Gerasimenko_Petrina}  for surveys, and  \cite{GSRT, Pulvirenti_Saffirio_Simonella, Pulvirenti_Simonella}.
In all these studies, the starting point is a system of evolution equations for the  correlation functions, which are finite dimensional projections of the probability distribution assigned on the whole particle system. 
The $k$-particle correlation function $F_k \big( t, (x_i,v_i)_{i \leq k} \big)$  describes 
the distribution at time $t$ of $k$ typical particles with positions denoted by $x_i$ and velocities by $v_i$.
These correlation functions obey the well known BBGKY hierarchy which states that, due to  binary collisions, the variation in time of $F_k$ depends on the distribution of $k+1$ particles, $F_{k+1}$. As a consequence, the density of a typical particle $F_1$ does not follow  a closed equation, and a substantial amount of work is required to prove that the correlation functions factorise in the Boltzmann-Grad limit, in such a way that the evolution equations can be closed with only the first limiting correlation function.

In the specific case of the billiard, the equations of the BBGKY hierarchy are singular and it is more convenient, mathematically, to represent the particle density in terms of an iterated Duhamel series (see e.g.\,\cite{Simonella_BBGKY, GSRT} on the justification of this series).  This formula relates the density of a typical particle at time $t$, $F_1(t)$, to the initial probability measure (described by the collection $\left(F_k(0)\right)_{k \geq 1}$),
by application of intertwined transport and collision operators. Ultimately this representation is studied by interpreting these chains of operators as integrals over  so-called ``pseudo-trajectories'' (\cite{GSRT}), which play the role of a characteristic flow.

In \cite{BGSS2}, the analysis of the correlation functions has been improved in order to control the fluctuations of the empirical measure and  not only its mean.
The key feature is the computation of the precise asymptotics for the exponential moment of the empirical measure, which is done via a cumulant expansion. From this, several results can be derived, namely the fluctuating Boltzmann equation and the large deviations for the hard sphere gas. We refer to \cite{BGSS1, BGSS_ICM}, and references therein, for a survey on these results and their physical interpretation, including the relation with stochastic particle dynamics (Kac process). 
One possible disadvantage of the iterated Duhamel representation is that the microscopic dynamics is not used directly, and  the above-mentioned pseudo-trajectories do not correspond to physical trajectories.
The link between the BBGKY hierarchy and physical trajectories is indeed very indirect (see for instance~\cite{BGS_bigparticle},  and~\cite{Simonella_BBGKY, Pulvirenti_Simonella_empirical}).
Furthermore, the pseudo-trajectories are time-oriented, and followed backwards up to the initial time, which may appear to be at variance with the naive idea of a stochastic process.

\medskip
The goal of the present  paper is  to obtain a statistical description of physical trajectories,   without  using the iterated Duhamel representation. More precisely, we will apply the cluster expansion method to real trajectories, controlling correlations of dynamical type. 
Cluster expansions have a long history in statistical mechanics, where they have been widely used to analyse correlations in Gibbs measures and  the equilibrium behaviour of particle systems; see e.g.\,\cite{Ruelle_livre, friedli2017statistical}. Originally, they were designed to study gases in the low density regime and to establish thermodynamic relations in the spirit of 
the virial expansion. This powerful method was then proved to be effective in more general contexts, whenever relevant observables can be identified as weakly (or rarely) interacting. For example, for the Ising model in the phase transition regime, the relevant observables are no longer the spins, but the spin contours which form a dilute gas of contours at low temperature (\cite{friedli2017statistical}).
More recently, the theory  was extended to an abstract framework covering continuous particle systems and polymers  \cite{poghosyan2009abstract}. Several surveys have been written on the many aspects of the cluster expansion (which we do not quote exhaustively) and the topic has been matter of investigation for decades \cite{Brydges_course, Kotecky_Preiss, Gallavotti_Bonetto_Gentile, Fernandez_Procacci, Winter_Simonella, Jansen_Kolesnikov}.

In our approach, the elementary observable is not just the position (and velocity) of a classical particle, but rather the trajectory (or path) in $[0,T]$ of a ``dynamical cluster'' (see Figure \ref{figure: cluster paths 3D} below). We show that for some short time $T>0$, the interactions are sufficiently rare so that the cluster expansion on such trajectories converges. Concerning correlations, the dynamical clusters studied in this paper play a similar role as the (perhaps more standard) cumulants, which were the main focus in \cite{BGSS2}: they grasp information on the dynamics on finer and finer scales. However, the cluster dynamics has an independent interest, when considering the problem of giving a mathematical meaning to the dynamics of the infinite system of particles (as obtained, for instance, from $N$ hard spheres in the Boltzmann-Grad limit). Indeed the cluster dynamics retains all the information on the limiting dynamics. The infinite system can be reinterpreted as a collection of clusters (groups of particles), moving independently and interacting randomly as in \cite{Sinai_cluster_1972, Sinai_cluster_1974}. The law of a typical cluster follows a coagulation process, the nonlinear behavior of which is reminiscent of Tanaka's process for the Boltzmann equation \cite{Tanaka} (see also \cite{Albeverio_Ruediger_Sunder_17} for space-dependent variants, and references therein). Our method provides a  rigorous derivation of the limiting cluster process,  Theorem \ref{Prop: convergence dynamical clusters}.

We mention that dynamical clusters have been investigated previously, by means of numerical experiments on $N$ hard spheres~\cite{Gabrielov_Keilis-Borok_Keilis-Borok_Zaliapin} (see also \cite{McFadden_Bouchard}). The main focus of these simulations is a phase transition in the cluster formation process: a giant (macroscopic) cluster appears abruptly at a critical time. Later on in \cite{Patterson_Simonella_Wagner1} it was shown that, in the Boltzmann-Grad limit, the cluster dynamics can be described by explicit formulas and that in this limit, the critical transition takes place in times of order one. Finally in \cite{Patterson_Simonella_Wagner2,heydecker2019bilinear}, the cluster process was studied for Kac particle systems, obtaining more precise results in spatially homogeneous cases.
Theorem \ref{Prop: convergence dynamical clusters} below derives rigorously the cluster process from hard spheres, giving in particular a proof to the formal statements in \cite{Patterson_Simonella_Wagner1}, at least for short times. 

\medskip
Note that, contrary to Lanford's approach  which studies correlations for finite dimensional projections, the cluster expansion is implemented here directly on real paths of the full microscopic system -- one will resort to projections only in a final step, in order to identify the Boltzmann equation for the limiting density (Section \ref{sec: Boltzmann equation}), or the asymptotics of the partition function governing the fluctuation theory and the large deviations (Section \ref{sec: Asymptotics free energy}). 
  Compared to \cite{BGSS2}, this provides therefore an alternative take on the Boltzmann-Grad limit, with twofold interest: (i) a more direct link with physical trajectories; (ii) a representation of observables in terms of a forward-in-time process, with randomness entering through the initial measure. We refer also to \cite{Matthies-Theil,Gerasimenko_Gapyak_21} for approaches sharing analogies with ours.

\bigskip

Let us now detail the model.  We consider a microscopic system of  identical hard spheres of unit mass and of diameter~$\eps$.
The  motion of $N$ such hard spheres is ruled by a system of ordinary differential equations, which are set in~$ ( \bbT ^d\times \bbR^d)^{N }$ where~$\mathbb T^d = [0,1]^d$ is  the unit~$d$-dimensional periodic box with $d \geq 2$: writing~${\textbf x}^{\e}_i \in  \bbT ^d$ for the position of the center of the particle labeled by~$i$ and~${\textbf v}^{\e}_i \in  \bbR ^d$ for its velocity, one has
\begin{equation}
\label{hardspheres}
{d{\textbf x}^{\e}_i\over dt} =  {\textbf v}^{\e}_i\,,\quad {d{\textbf v}^{\e}_i\over dt} =0 \quad \hbox{ as long as \ } |{\textbf x}^{\e}_i(t)-{\textbf x}^{\e}_j(t)|>\eps  
\quad \hbox{for \ } 1 \leq i \neq j \leq N
\, ,
\end{equation}
with specular reflection at collisions, i.e. when $|{\textbf x}^{\e}_i(t)-{\textbf x}^{\e}_j(t)|=\eps$ 
\begin{equation}
\label{defZ'nij}
\begin{cases}
 \left({\textbf v}^{\e}_i\right)'
:= {\textbf v}^{\e}_i -  \big( ({\textbf v}^{\e}_i - {\textbf v}^{\e}_j)\cdot \pmb{\omega} \big) \, \pmb{\omega}   \\
\left({\textbf v}^{\e}_j\right)' := {\textbf v}^{\e}_j + \big( ({\textbf v}^{\e}_i - {\textbf v}^{\e}_j)\cdot \pmb{\omega} \big) \, \pmb{\omega} 
\end{cases}
\end{equation}
where $\pmb{\omega} := ({\textbf x}^{\e}_i (t) - {\textbf x}^{\e}_j (t))/\eps$ is the unit vector pointing along the relative position at the collision time $t$. 
The collections of~$N$ positions and velocities are denoted respectively by~$X_N := (x_1,\dots,x_N) $ in~$ \bbT^{dN}$ and~$V_N := (v_1,\dots,v_N) $ in~$ \bbR^{d N}$,  and we set~$Z_N:= (X_N,V_N) $, with~$Z_N =(z_1,\dots,z_N)$, $z_i = (x_i,v_i)$. 
A set of~$N$ particles is characterised by the time-zero configuration~${\textbf Z}^{\eps 0}_N =  ({\textbf z}^{\eps 0}_1,\dots , {\textbf z}^{\eps 0}_N)$   in the phase space 
\begin{equation}
\label{eq: D-def}
{\mathcal D}^{\eps}_{N} := \big\{Z_N \in (\bbT ^d\times\bbR^d)^N \, / \,  \forall i \neq j \, ,\quad |x_i - x_j| > \eps \big\} \, ,
\end{equation}
and an evolution
$$ 
t \mapsto {\textbf Z}^{\eps}_N(t) = \big({\textbf z}^{\eps}_1(t),\dots , {\textbf z}^{\eps}_N(t)\big)\;,\qquad t>0
$$
according to the  flow \eqref{hardspheres}-\eqref{defZ'nij} (well defined 
almost surely under the Lebesgue measure \cite{Alexander}).

\medskip

The dynamics is deterministic, but the initial data are chosen randomly according to the grand canonical formalism described below (see \cite{Ruelle_livre} for details).
The number of particles $\cN$ is a random variable so that 
the initial measure is defined  on the phase space 
$$
\cD^\eps := \bigcup_{N \geq 0} \cD^\eps_N
$$
(notice that $\cD^\eps_N = \emptyset$ for $N$ large due to the exclusion condition). 
Initially,  the probability density of finding $N$ particles with configuration $Z_N$ is given by
\begin{equation}
\label{eq: initial measure}
\frac{1}{N!}
W^{\eps }_{N}(Z_N) 
:= \frac{1}{\cZ^ \eps} \,\frac{\mu_\eps^N}{N!} \, 
\indc_{{\mathcal D}^{\eps}_{N}} (Z_N) \; \prod_{i=1}^N f^0 (z_i) \, , 
\end{equation} 
where the distribution of a single particle $f^0$ is a Lipschitz continuous probability density  on $\bbT ^d\times\bbR^d$ satisfying the following bound
for some constants $\beta, C_0 >0$
\begin{equation}
\label{eq: Gaussian bound}
\forall z \in \bbT ^d\times\bbR^d, \quad 
f^0(z) \leq C_0 \cM_\beta (v)
\quad \text{with} \quad 
\cM_\beta (v) := \frac{\beta^{d/2}}{(2\pi)^\frac d2} \exp \left( - \beta \frac{ |v|^2}{2} \right),
\end{equation}
and the partition function is given by
\begin{equation}
\label{eq: partition function}
\cZ^\eps :=  1 + \sum_{N\geq 1}\frac{\mu_\eps^N}{N!}  
\int_{\cD^\eps_N} d Z_N \prod_{i=1}^N  f^0(z_i)  \,.
\end{equation}
The stated assumptions on the initial measure will be kept throughout all the paper.
The probability of an event~$A$ with respect to the measure \eqref{eq: initial measure} 
will be denoted~${\mathbb P}_\eps(A)$, and~${\mathbb E}_\eps$ will be the expectation.

In the low density regime, referred to as  the Boltzmann-Grad scaling, the expectation of the number of particles  $\cN$ is tuned by the parameter 
\begin{equation}
\label{eq: Boltzmann-Grad}
\mu_\e := \e^{-(d-1)},
\end{equation} ensuring that the mean free path between collisions is of order one \cite{Grad49}. 
Definition \eqref{eq: initial measure} implies that 
$$
\lim_{\eps \to 0} \; \mu_\e^{-1}{\mathbb E}_\eps\left(\cN\right) =  1.
$$ 
Let us denote by~$W^{\eps }_{N}(t,Z_N) $ the    probability density of finding $N$ particles with configuration $Z_N$  at time~$t$, which satisfies the Liouville equation
\begin{equation}
\label{Liouville}
	\partial_t W^{\eps}_N +V_N \cdot \nabla_{X_N} W^{\eps}_N =0  \,\,\,\,\,\,\,\,\, \hbox{on } \,\,\,{\mathcal D}^{\eps}_{N}\, ,
\end{equation}
	with specular reflection on the boundary. The  $k$-th correlation function is defined by
	\begin{align*}
F_k^{\eps} (t,Z_k)
:= \mu_\eps^{-k} \,
\sum_{p=0}^{\infty} \,\frac{1}{p!}\, \int dz_{k+1}\dots dz_{k+p} \,
W_{k+p}^{\eps } (t,Z_{k+p}) 
\end{align*}
and as
mentioned earlier in the introduction, the key result originally derived by Lanford \cite{Lanford} is the convergence of~$F_k^{\eps}$
 to the tensor product $ f^{\otimes k}$ where $f$ is the solution of the Boltzmann equation with initial data $f^0$
\begin{equation}
\left\{ \begin{aligned}
& \partial_t f +v \cdot \nabla _x f 
= \! \displaystyle\int_{\bbR^d}\int_{{\mathbb S}^{d-1}}   \! \Big(  f(t,x,w') f(t,x,v') - f(t,x,w) f(t,x,v)\Big) 
\big ((v-w)\cdot \omega\big)_+ \, d\omega \,dw \,  ,\\
&  f(0,x,v) = f^0(x,v)
 \end{aligned}
  \right. 
\label{eq:Beq}
  \end{equation}
with   precollisional velocities $(v',w')$ defined by the scattering law
\begin{equation}
\label{eq: scattlaw}
 v' := v- \big( (v-w) \cdot \omega\big)\,  \omega \,  ,\qquad
 w' :=w+\big((v-w) \cdot \omega\big)\,  \omega  \, .
\end{equation}
This can be rephrased in terms of the convergence of  the  empirical density
\begin{equation}
\label{eq: empirical}
\pi^{\eps}_t :=\frac{1}{\mu_\e}\sum_{i=1}^\cN \delta_{{\textbf z}^\e_i(t)}\;.
\end{equation} 
The convergence of the particle system to the Boltzmann equation can indeed be understood as follows.
\begin{thm}[Lanford, \cite{Lanford}]
\label{thm: Lanford}
There exists a time $T_L >0$ such that 
for any test function $h: \bbT ^d\times \bbR^d \to \bbR$, any $\delta>0$ and  $t \in [0,T_L]$,
\begin{equation}
\label{eq: convergence Boltzmann}
\bbP_\eps \left( \Big|\pi^\eps_t(h) -  
   \int_{\bbT ^d\times \bbR^d} dz f(t,z) h(z)\Big| > \delta \right) \xrightarrow[\mu_\eps \to \infty]{} 0 \,.
\end{equation}
The time $T_L$ depends only on the smooth  function $f^0$ through the constants~$C_0$ and~$\beta$ appearing in~\eqref{eq: Gaussian bound}.
\end{thm}
Notice that stronger convergence statements can be found \cite{IP89,Cercignani_Gerasimenko_Petrina,Cercignani_Illner_Pulvirenti,GSRT, Pulvirenti_Simonella, denlinger, BGSS,Gerasimenko_Gapyak_18,Gerasimenko_Gapyak_20}.
All these studies rely on the BBGKY hierarchy and one of the goals of this paper is to provide 
an alternative derivation of Theorem \ref{thm: Lanford} by applying directly a cluster expansion at the level of the particle system.
More generally, we are interested in the whole path of particle trajectories during a given time interval $[0,T]$.
Let $\bbD ([0,T], \bbT^d \times \bbR^d)$ be the set of single particle trajectories ${\textbf z}^\e ([0,T])$, which are functions piecewise linear continuous in position and piecewise constant in velocity.
Then the generalised empirical measure is defined by
\begin{equation}
\label{eq: empirical general}
\pi^{\eps}_{[0,T]} :=\frac{1}{\mu_\e}\sum_{i=1}^\cN \delta_{{\textbf z}^\e_i([0,T])} \;.
\end{equation} 
To derive sharp estimates  on the  empirical measure, it is key to control its exponential moments, i.e.\,the Laplace transform  
\begin{equation}
\label{eq: log partition function}
\widetilde \gL^\eps_T (e^h) 
:= \frac{1}{\mu_\eps}\log \bbE_\eps \left[ \exp \Big(   \mu_\eps\,\pi^\eps_{[0,T]}(h)  \Big)   \right] 
= \frac{1}{\mu_\eps}\log \bbE_\eps \Big[ \exp \big(   \sum_{i=1}^{\cN} h\big( {\textbf z}^{\eps}_i([0,T]) \big)   \big)\Big]  \,,
\end{equation}
for test functions  $h : \bbD ([0,T], \bbT^d \times \bbR^d) \to \bbC$ measuring information on a single particle trajectory ${\textbf z}^\eps ([0,T])$. 
Indeed, it is well known in probability theory that the large deviations can be related to the Legendre transform of 
$h \mapsto \widetilde \gL^\eps_T (e^h)$ and that the fluctuations are coded by the characteristic function which amounts to considering functions $h$  of size~$O(\mu_\eps^{-\frac12})$ with imaginary values. We refer to \cite{BGSS2, BGSS1} for the derivation of the fluctuating Boltzmann equation and of the large deviations once the asymptotic behaviour of $\widetilde\gL^\eps_T$ has been characterised.

In the rest of this paper, we first analyse, in Section \ref{sec: Dynamical cluster expansion}, several types of dynamical interactions and show that the static cluster expansion on equilibrium measures can be extended to a dynamical cluster expansion in the dilute regime. Then we study its Boltzmann-Grad limit in Section \ref{sec: Boltzmann-Grad limit}. Section \ref{sec: limiting equations} is devoted to the derivation of limiting dynamical equations, including the Boltzmann equation as in Theorem \ref{thm: Lanford}, and a coagulation-type equation driving the evolution of the dynamical clusters (Theorem \ref{Prop: convergence dynamical clusters}).


\medskip

\section{Dynamical cluster expansion}
\label{sec: Dynamical cluster expansion}

\subsection{Decomposition in cluster paths} 
\label{trajectory-subsec}
 
Throughout this section, we study the hard sphere dynamics on a fixed time interval $[0,T]$
and implement a cluster expansion to study the functional \eqref{eq: log partition function}
 as well as a natural extension which will be introduced in 
\eqref{eq: log modified partition function}.  This corresponds to a rather standard statistical mechanics computation, if not for the fact that the positions of classical particles are now replaced by their paths.
Even though the gas is extremely dilute in the Boltzmann-Grad asymptotics, particles are likely to interact dynamically.  In this case their trajectories are strongly modified by scattering.
Thus  to implement the dynamical cluster expansion, a good point of view is to first group particles which undergo collisions during $[0,T]$ and then to perform the standard cluster expansion on these groups of particles.

\smallskip

\begin{figure}[h] 
\centering
\includegraphics[width=3.2in]{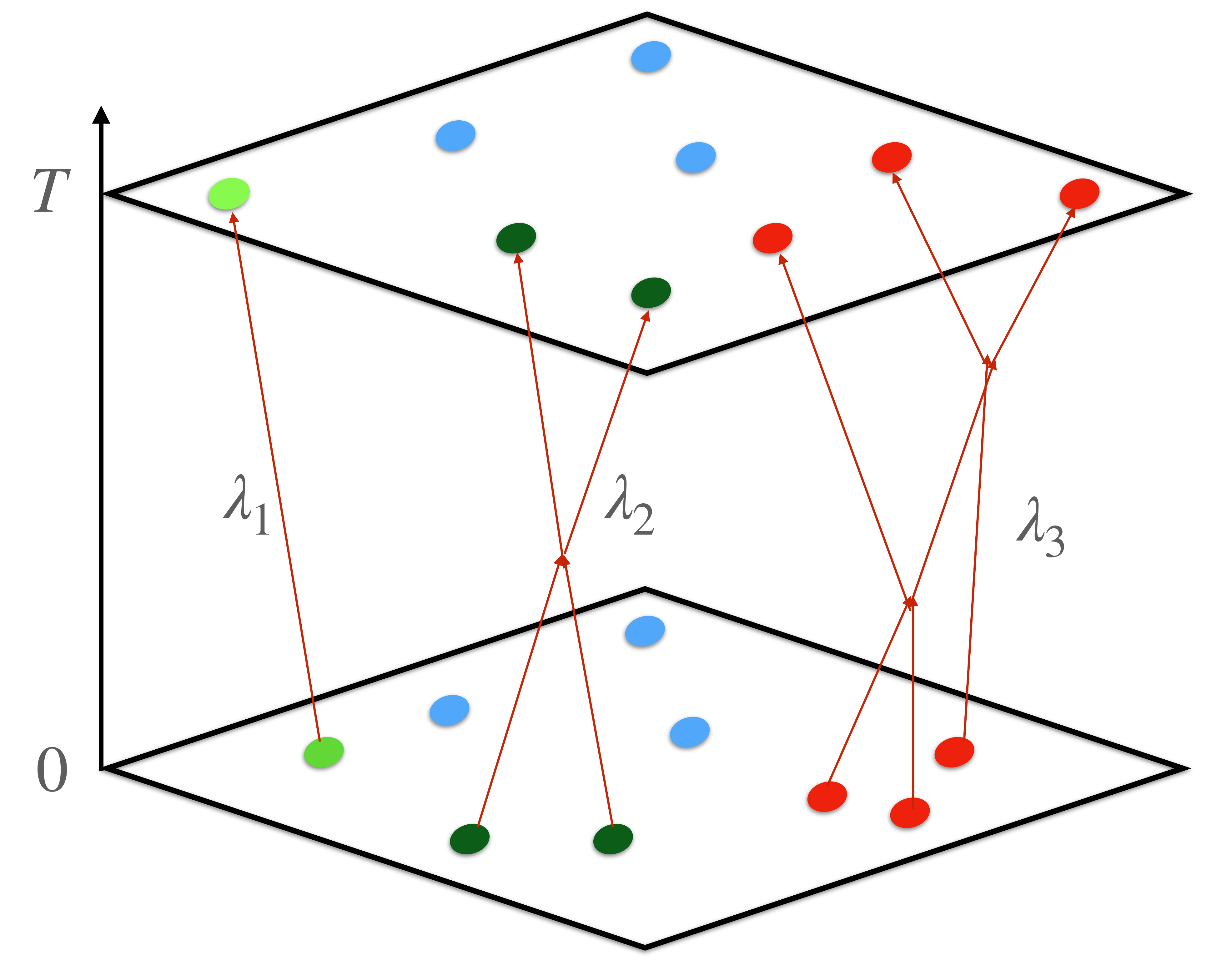} 
\caption{\small 
On this figure,  three $\eps$-cluster paths $\lambda_1, \lambda_2, \lambda_3$ are depicted in the time interval $[0,T]$. They are formed by groups of particles interacting dynamically.
As the $\eps$-cluster paths do not intersect, the trajectories within an $\eps$-cluster path are only determined by the hard sphere dynamics restricted to the particles in this cluster path.
Note that the other blue particles will form more cluster paths which have not been represented.} \label{figure: cluster paths 3D}
\end{figure}

\begin{Def}[$\eps$-cluster path]
\label{def: cluster paths}
We say that two particles interact dynamically on  $[0,T]$ if they collide on that time interval. Given a set of particle trajectories, a graph of dynamical interactions is built by adding an edge~$\{i,j\}$ if two particles~$i$ and~$j$ from that set interact dynamically on  $[0,T]$.

 An  \emph{$\eps$-cluster path} $\lambda $ on $[0,T] $  is   a set~$ \{ {\textbf z}^\eps_i ([0,T]) \}_{i \in I_{\lambda }}$ of  particle trajectories having a connected graph of dynamical interactions, and which do not interact dynamically with   particles outside that set.
 The number of particles in the $\eps$-cluster path  $\gl  $ is denoted by~$| \gl  |:=|I_{\lambda }|$.
 \end{Def}
 
\begin{Rmk}\label{rmk parametrization}As the hard sphere dynamics is deterministic, the dynamical condition to form the $\eps$-cluster path $\gl $ on $[0,T] $  is coded in the initial data.
\end{Rmk}

By definition, particles in an $\eps$-cluster path $\gl$ never collide with particles which do not belong to $\gl$. 
Thus if a trajectory of the whole microscopic system~${\textbf Z}^\eps_N ([0,T])$ is decomposed into a  partition $\{ \gl_1, \dots, \gl_k \}$ of $\eps$-cluster paths, then particles in different cluster paths $\gl_i,\gl_j$ do not collide. This condition is denoted by $\gl_i \not \sim_\eps  \gl_j$. Given the time interval $[0,T]$, 
the decomposition of~${\textbf Z}^\eps_N ([0,T])$ into $\eps$-cluster paths is unique (see Figure \ref{figure: cluster paths 3D}).
Note that  the total kinetic energy of an $\eps$-cluster path $\gl_j$ is  time independent within the time of existence of the path: $\displaystyle {\mathcal E}_{\gl_j}:=\frac12\sum_{i \in I_{\gl_j}} |v^\eps_i(t)|^2$ does not depend on~$t \in [0,T]$.

By  Definition \ref{def: cluster paths},  given a time $T>0$, one has the following partition of unity
\begin{equation}
\label{eq:cluster path decomposition}
 \indc_{ \mathcal D^\eps_N }  =  \sum_{k \leq N} \sum_{\{ \gl_1, \dots, \gl_k\} \in \cP_{{\textbf Z}^\eps_N}^k} \left( \prod_{j \neq j'} \indc _{ \lambda_j \not\sim_\eps \lambda_{j'}} \right) \left( \prod_{j= 1} ^{k}  \indc_{ {\lambda_j} \hbox{ \scriptsize$\eps$-cluster path } } \right),
\end{equation}
where $\cP_{{\textbf Z}^\eps_N}^k$ is the set of partitions of the single particle trajectories~$({\textbf z}^\eps_i ([0,T]))_{i\leq N} $ into $k$ sets.  We stress that, in formula \eqref{eq:cluster path decomposition}, the $\eps$-cluster paths are seen as functions of the initial (time-zero) configuration $\left( z_i\right)_{i \in I_{\lambda_j}}$.

 This formula with the dynamical exclusion $ \prod_{j \neq j'} \indc _{ \lambda_j \not\sim_\eps \lambda_{j'}} $ suggests to consider the system as a gas of $\eps$-cluster paths. We expect this gas to be in   a dilute regime for small times (roughly of the order of the convergence time $T_L$ introduced in Theorem~\ref{thm: Lanford}). The main results of this paper will actually  be obtained by applying classical cluster expansions to this gas of $\eps$-cluster paths.
 
\bigskip

We indicate by~${\textbf Z}:={\textbf Z}([0,T])$ a generic trajectory of size~$|{\textbf Z}|:=n$,  namely any element of~${\mathbb D}^n([0,T],{\mathbb T}^d\times {\mathbb R}^d)$. Dropping the dependence on~$T$ when there is no ambiguity, we have that ${\textbf Z} = \left\{ {\textbf z}_i  ( [0,T]  ) \right\}_{i \leq n}$, and we assume that the  trajectory conserves the total kinetic energy, meaning~$\cE_{\textbf Z}\displaystyle  :=\frac12\sum_{i\leq n} |{\textbf v}_i(t) |^2$ does not depend on~$t\in [0,T]$. We equip the set of trajectories with the following norm: a sequence of trajectories~${\textbf Z}^{\ell}$ converges to a trajectory~${\textbf Z}   $ as~$\ell$ goes to infinity if~$ | {\textbf Z}^\ell |= |{\textbf Z}|$ for~$\ell$  large and
\begin{equation}
\label{eq: norm T}
\|{\textbf Z}^\ell - {\textbf Z}\|_T :=  \sup_{t\in[0,T]} \, \sup_{ i\leq  |{\textbf Z}|} \big| {\textbf z} ^\ell_i(t) - {\textbf z}_i(t) \big| 
\xrightarrow[\ell \to \infty]{} 0\;.
\end{equation}

\bigskip

Let us now consider a complex-valued test function $H$ defined on  the set of {trajectories}, 
such that the following bound holds uniformly over trajectories 
 conserving the total kinetic energy
\begin{equation}
\label{eq: hypotheses h}
\left|  H({\textbf Z})  \right| \leq   c_1|{\textbf Z}| +c_2 \cE_{\textbf Z}
\quad \text{for some fixed constants $c_1>0$,  \quad $0<c_2 \leq \beta/4$}
\,,
\end{equation}
where we recall that $|{\textbf Z}|=n$ indicates the total number of particles involved in the trajectory.
To study the Boltzmann-Grad limit in Section \ref{sec: Boltzmann-Grad limit}, 
we will need  $H$ to be continuous for the   uniform convergence notion \eqref{eq: norm T}: if~${\textbf Z}^\ell$ converges to~${\textbf Z}   $ as~$\ell$ goes to infinity, then
\begin{equation}
\label{eq: hypotheses h continuity}
\lim_{\ell \to \infty} H( {\textbf Z}^\ell) = H({\textbf Z}) \, .
\end{equation}
We shall obtain information about the empirical measure $\pi^\eps_{[0,T]}$ defined in (\ref{eq: empirical general}) by choosing~$H$ of the form
\begin{equation}
\label{eq: hypotheses h bis}
H({\textbf Z}) = \sum_{i =1}^{|{\textbf Z}|} h   \big( {\textbf z}_i ([0,T]) \big) \,,
\end{equation}
but we can actually consider much more general functions $H$ satisfying (\ref{eq: hypotheses h}), which will give access to microscopic events such as   the size or  the structure of cluster paths.
In fact, 
the empirical measure $\pi^{\eps}_{[0,T]}$ introduced in \eqref{eq: empirical general} can be extended to the space of $\eps$-cluster paths 
\begin{equation}
\label{eq: modified empirical general}
\Pi^{\eps}_{[0,T]} :=\frac{1}{\mu_\e}\sum_i \delta_{ \gl_i} \;,
\end{equation} 
 as well as the exponential moments   
\begin{equation}
\label{eq: log modified partition function}
\gL^\eps_T (e^H) 
:= \frac{1}{\mu_\eps}\log \bbE_\eps \left[ \exp \Big(   \mu_\eps\, \Pi^\eps_{[0,T]}(H)  \Big)   \right] 
= \frac{1}{\mu_\eps}\log \bbE_\eps \Big[ \exp \big(   \sum_i H \big( \gl_i  \big)   \big)\Big]  \,.
\end{equation}
Notice that~(\ref{eq: log partition function}) is a specific case of~(\ref{eq: log modified partition function}) for~$H$ of the form~(\ref{eq: hypotheses h bis}).

\medskip

  Recalling the definition~(\ref{eq: partition function}) of the partition function and the decomposition~(\ref{eq:cluster path decomposition}) into $\eps$-cluster paths, we define
 the modified partition function
\begin{equation}
\label{eq: modified partition function}
\cZ^\eps_T (e^H)  :=  1 + \sum_{N\geq 1}\frac{\mu_\eps^N}{N!}  
\sum_{k \leq N} 
\int_{\bbT^{d N} \times \bbR^{dN}} d Z_N  \;
\sum_{\{  \gl_1, \dots, \gl_k\} \in \cP_{{\textbf Z}^\eps_N }^k} \left(  \prod_{j \neq j'} \indc _{ \lambda_j \not\sim_\eps \lambda_{j'}}  \right) 
 \left( \prod_{j= 1}^k  \indc_{ {\lambda_j} \hbox{ \scriptsize$\eps$-cluster path } } F ^{(H)} (\gl_j) \right) \,  , 
\end{equation}
where the weight of  a generic trajectory~$ {\textbf Z}$
is defined by  
\begin{equation}
\label{eq: weight cluster path}
F ^{(H)} ({\textbf Z}) := e^{H({\textbf Z}) }\prod_{i \leq |{\textbf Z}|} f^0({\textbf z}_i(0)), 
\end{equation}
  and we simply write $F ^{(H)} (\gl)$ when the trajectory is an~$\eps$-cluster path $\gl$.

We stress the fact that, by definition of the $\eps$-cluster paths $\lambda_{ j}$ appearing in~(\ref{eq: modified partition function}), 
the particles in $\lambda_{ j}$ are not allowed to overlap initially  
$$
\forall i \not =  i' \in \{1,\dots,|\lambda_{ j}|\} \, , \qquad  |x_i - x_{i'}| \geq \eps,
$$
otherwise the particle trajectories would be ill defined.
The exponential moment of the $\eps$-cluster paths can be rewritten as  
\begin{equation}
\label{eq: log modified partition function bis}
\gL^\eps_T (e^H) 
 = \frac{1}{\mu_\eps} \left( \log \cZ^\eps_T (e^H) - \log \cZ^\eps_T (1) \right)  \,.
\end{equation}

\begin{Rmk}
For $H$ of the form \eqref{eq: hypotheses h bis}, the exponential moment \eqref{eq: log partition function}
can be recovered $\gL^\eps_T (e^H) =  \widetilde \gL^\eps_T (e^h)$.
It can be reformulated in terms of  the modified partition function
\begin{equation}
\label{eq: modified partition function bis}
\widetilde \cZ_T^\eps (e^h) :=  1 + \sum_{N\geq 1}\frac{\mu_\eps^N}{N!}  
\int_{\cD^\eps_N} d Z_N  \prod_{i=1}^N  f^0(z_i) \exp \Big(  h ( {\textbf z}^{\eps}_i ([0,T]) \Big)  
\end{equation}
so that 
\begin{equation}
\label{eq: log partition function bis}
\widetilde \gL^\eps_T (e^h) 
 = \frac{1}{\mu_\eps} \left( \log \widetilde\cZ^\eps_T (e^h) - \log \widetilde\cZ^\eps_T (1) \right)  \,.
\end{equation}

\end{Rmk}

\bigskip

 In~(\ref{eq: modified partition function})
the integration variables are the initial particle configurations $Z_N$ (which as recalled in Remark~\ref{rmk parametrization} fix the dynamics on $[0,T]$ and therefore the partition into cluster paths). Fubini's theorem enables us to consider now as the variables of interest   the cluster paths $\gl_j$. 
It is convenient to define the following integration measure~:
\begin{equation}
\label{eq: distribution cluster path}
d \mes ( \gl ) := \frac{\mu_\eps^{| \gl |}}{| \gl | \, !}\,  \indc_{ {\lambda} \text{ \scriptsize $\eps$-cluster path on $[0,T]$} }\, F^{(H)} (\gl)
\;  d Z_{|\lambda|} \, , 
\end{equation}
where the support of the measure is restricted by the dynamical constraint $\indc_{ {\lambda} \text{ $\eps$-cluster path on $[0,T]$} }$, meaning 
that the trajectories associated with the initial data~$Z_{|\lambda|}$ form an $\eps$-cluster path in the time interval $[0,T]$.
 Later on, the time $T$ will be chosen small enough so that the gas of cluster paths is dilute 
and the cluster expansion converges. 
We also define~$d |\mes|$ as in \eqref{eq: distribution cluster path} with the modulus $| F^{(H)} |$ to take into account test functions $H$ with complex values.
With these notations,  the partition function \eqref{eq: modified partition function} can be rewritten as a partition function of a gas of cluster paths interacting by exclusion.
Indeed  a particle configuration 
$Z_N$ is partitioned into $k$ cluster paths with cardinalities $n_1, \dots , n_k$
$$
\cZ^\eps_T (e^H)  =  1 + \sum_{N\geq 1}
\frac{1}{N!}  
\sum_{k \leq N}  \frac{1}{k!}
\sum_{n_1, \dots , n_k \atop n_1 + \dots + n_k = N}
\frac{N!}{n_1! \dots n_k!}  \mu_\eps^{n_1 + \dots + n_k}
\int_{\bbT^{d N} \times \bbR^{dN}} \prod_{j= 1}^k \left(  d  Z_{|\lambda_j|}
  \indc_{ {\lambda_j} \text{ \scriptsize$\eps$-cluster path on $[0,T]$} \atop \text{\scriptsize of size $n_j$}} F^{(H)} (\gl_j) \right) \; \prod_{j \neq j'} \indc _{ \lambda_j \not\sim_\eps \lambda_{j'}} 
$$
where $\frac{N!}{n_1! \dots n_k!}  $ is the number of partitions of~$\{1,\dots,N\}$ into~$k$ ordered components of cardinalities $n_1, \dots ,n_k$, and 
$\frac{1}{k!}$ removes the multiple counting due to the order.  
By Fubini's theorem, we then  find
\begin{equation}
\cZ^\eps _T(e^H)  
 =  1 + \sum_{k \geq 1}   \frac{1}{k!} \int   
\prod_{j= 1}^k \sum_{n_j} 
\;  \left(  \frac{\mu_\eps^{n_j}}{ n_j!}   d Z_{|\lambda_j|}
  \indc_{ {\lambda_j} \text{ \scriptsize$\eps$-cluster path on $[0,T]$} \atop \text{\scriptsize of size $n_j$}} F^{(H)} (\gl_j) \right) \; \prod_{j \neq j'} \indc _{ \lambda_j \not\sim_\eps  \lambda_{j'}}\;.
  \label{eq: partition cluster path dvp}
\end{equation}
Note that the sums can be exchanged since they are finite  (the total number of particles in the box is finite for any fixed $\eps>0$ thanks to the exclusion condition). 

Finally, with notation \eqref{eq: distribution cluster path}, we can write
\begin{equation}
\label{eq: partition gas cluster paths}
\cZ^\eps_T (e^H) =  1 + \sum_{k \geq 1}\frac{1}{k!}  
\int d \mes ( \gl_1) \dots d \mes ( \gl_k ) \; 
\prod_{j \neq j'} \indc _{ \lambda_j \not\sim_\eps \lambda_{j'}}
\end{equation}
where the integration for each $\eps$-cluster path $\lambda$ is with respect to the number~$|\lambda|$ of particles in the $\eps$-cluster path as well as the initial coordinates~$ Z_{|\lambda|}$  of the particles (see \eqref{eq: partition cluster path dvp}).

\subsection{Cluster  expansion on  $\eps$-cluster paths} 
\label{cluster-subsec}

In the theory of cluster expansions, it is customary to control 
the exclusion interaction in \eqref{eq: partition gas cluster paths} 
$$
\prod_{j \neq j'} \indc _{ \lambda_j \not\sim_\eps \lambda_{j'}} = 
\prod_{j \neq j'}  \Big( 1 - \indc _{ \lambda_j \sim_\eps \lambda_{j'}} \Big)
$$ 
by expanding the product above and then rearranging the terms.
This requires to plug the $\eps$-cluster paths on an extended space, introducing a new type of dynamical interaction (called \emph{overlap}) and further combinatorial decompositions. Let $\pmb{\cG}_k$ be the set of graphs with $k$ vertices and $\pmb{\cC}_k \subset \pmb{\cG}_k$ the subset of connected graphs. We can encode the exclusion between cluster paths by graphs~:
\begin{equation}\label{graph-dec}
\prod_{j \neq j'} \indc _{ \lambda_j \not\sim _\eps \lambda_{j'}} = 
\sum_{\cG \in \pmb{\cG}_k} \prod _{\{ j,j'\}  \in E (\cG)} (-\indc_{\lambda_{ j} \sim_\eps \lambda_{{j'}}})\,.
\end{equation}
This formula is made precise by the following definition.

\begin{Def}[Overlap and $\eps$-aggregate]
\label{def:overlap}
Consider a group of~$k$ $\eps$-cluster paths~$\{\lambda_1,\dots,\lambda_k\}$.    We say that two  $\eps$-cluster paths~$\lambda_j, \lambda_{j'}$  \emph{overlap}  if  two particles 
from  $\lambda_j$ and~$\lambda_{j'}$ are at a distance less than~$\eps$ at some time in $[0,T]$. We   write $\lambda_j \sim_\eps \lambda_{j'}$.

The $\eps$-cluster paths $\{ \gl_1, \dots, \gl_k \}$ form an \emph{$\eps$-aggregate} if the dynamical
graph with~$k$ vertices,  built by adding an edge each time two $\eps$-cluster paths overlap, is connected.
 The combinatorial function associated with this $\eps$-aggregate of size  $k$  is defined as 
\begin{equation}
\label{eq: phi}
\varphi_\eps (\gl_1, \dots, \gl_k) = 
\begin{cases}
1, & \text{if  \quad $k =1$}\, ,\\
\sum_{\cG \in \pmb{\cC}_k} \prod_{\{ j,j'\}  \in E (\cG)} (-\indc_{\lambda_ j \sim_\eps \lambda_{j'} }),  
 & \text{if  \quad  $k \geq 2$}\, ,
\end{cases}
\end{equation}
where the vertices of the graph $\cG$ are indexed by $\{ \gl_1, \dots, \gl_k \}$ and 
the product is over all the edges $E (\cG)$ of $\cG$. 
\end{Def}

\begin{figure}[h] 
\centering
\includegraphics[width=3.2in]{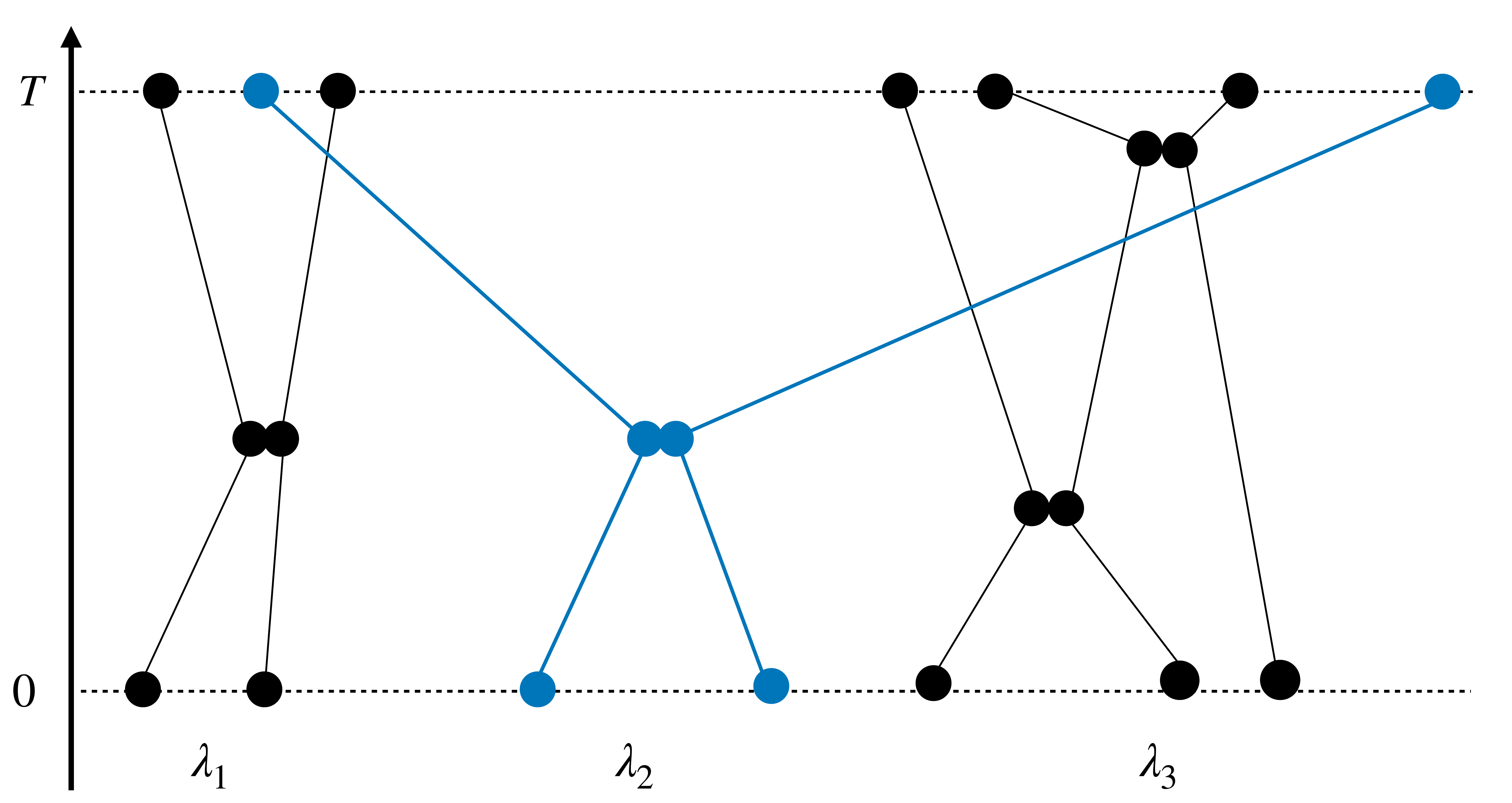} 
\caption{\small 
On this figure,  three  $\eps$-cluster paths~$\lambda_1, \lambda_2, \lambda_3$ forming an~$\eps$-aggregate are depicted (the $\eps$-cluster path $\lambda_2$ is drawn in blue only for the sake of clarity). 
The particle dynamics inside each cluster path is not modified by an overlap. 
Notice that the $\eps$-cluster paths $\lambda_2, \lambda_3$ overlap three times, this will be interpreted as a multiple edge in the generalized dynamical interaction graph introduced in Section \ref{sec: Boltzmann-Grad limit}.
}\label{figure: overlaps}
\end{figure}

\begin{Rmk}We stress the fact that an $\eps$-overlap between two $\eps$-cluster paths does not change the dynamics of the $\eps$-cluster paths : the particle trajectories remain encoded only by the particles within the cluster paths (see Figure \ref{figure: overlaps}). In this sense, an~$\eps$-overlap is a mathematical artefact which cannot be observed physically.
\end{Rmk}

Within the previous framework, the cluster expansion theory (see e.g. \cite{ueltschi2003cluster,poghosyan2009abstract})
leads to the following result.
\begin{Prop}
\label{prop: cluster expansion}
There exists a time $T>0$ and a constant $C$  depending only on the initial data~$f^0$ and on the constants~$c_1, c_2$ appearing in~(\ref{eq: hypotheses h}) such that  uniformly in $\eps$ small enough 
and for $k \geq 1$
\begin{equation}
\label{eq: bornes cluster expansion}
\int d | \mes| ( \gl_1 ) \dots d |\mes| ( \gl_k ) \; |\varphi_\eps (\gl_1, \dots, \gl_k )  |
\leq 
 \mu_\eps \; k! \; C^k \big(  T+ \eps \big)^{k-1}
\end{equation} 
and the cluster expansion of the modified partition function converges 
\begin{equation}
\label{eq: cluster expansion}
\log \cZ^\eps_T (e^H) =   \sum_{k \geq 1}\frac{1}{k!}  
\int d \mes ( \gl_1 ) \dots d \mes ( \gl_k ) \; \varphi_\eps (\gl_1, \dots, \gl_k )   \, .
\end{equation}
 \end{Prop}

\begin{Rmk}
The time $T>0$ for the convergence of the cluster expansion is of the same 
order as the convergence time $T_L$ to the Boltzmann equation in Theorem \ref{thm: Lanford}.
In both cases, it is a fraction of the mean collision time between two particles.
Indeed the requirement that the gas of $\eps$-cluster paths is dilute means that the  $\eps$-cluster path sizes
have to remain small.
A crude analogy can be made with an  Erd\H{o}s-Renyi graph built by choosing randomly edges 
among $N$ points with probability $T/N$.  For $T<1$, this procedure leads with high probability 
to a collection of small graphs which corresponds to the dilute phase we have in mind for the 
hard sphere dynamics. Instead, as soon as $T>1$, a macroscopic connected graph appears. We will make this analogy more precise in Section \ref{sec: tanaka} when analysing the dynamical clustering process.
We refer to \cite{Gabrielov_Keilis-Borok_Keilis-Borok_Zaliapin, McFadden_Bouchard, Patterson_Simonella_Wagner1, Patterson_Simonella_Wagner2, heydecker2019bilinear} for previous works 
on this dynamical phase transition in the case of particle dynamics. 
Reaching longer time asymptotics requires therefore new ideas and techniques.  
For this reason, we have made no attempt to optimise the time convergence $T>0$ in Proposition \ref{prop: cluster expansion}. 
\end{Rmk}

\medskip

\noindent
{\it Proof of Proposition \ref{prop: cluster expansion}.}

Assuming the validity of \eqref{eq: bornes cluster expansion}, the cluster expansion \eqref{eq: cluster expansion} follows from \cite{ueltschi2003cluster}. 
We sketch the proof of \eqref{eq: cluster expansion} below for the sake of completness.
Expanding the exclusion with (\ref{graph-dec})  in \eqref{eq: partition gas cluster paths}, we get
$$
\begin{aligned}
\cZ^\eps_T (e^H) &=   1+ \sum_{k \geq 1}\frac{1}{k!}  
\int d \mes ( \gl_1 ) \dots d \mes ( \gl_k ) \; 
 \sum_{\cG \in \pmb{\cG}_k} \prod_{ \{j,j'\} \in E( \cG) }   
\big(- \indc _{ \lambda_j \sim_\eps \lambda_{j'}}  \big)  \, .
\end{aligned}
$$
Any graph $\cG \in \pmb{\cG}_k$ can be decomposed 
into connected graphs~$\cG= \{ \cG_1, \dots, \cG_n\}$
with $| \cG_\ell| = m_\ell$ and $m_1 + \dots + m_n = k$. To do this, we partition $\{1, \dots , k\}$
into $n$ sets and then enumerate the graphs on each set
\begin{align*}
\cZ^\eps _T(e^H)
&= 1+  \sum_{k = 1}^\infty \frac{1}{k !}  
\sum_{n = 1}^k \frac{1}{n!}
\sum_{m_1, \dots , m_n \geq 1 \atop   m_1 + \dots + m_n = k} \frac{k!}{m_1 ! \dots m_n !}
\prod_{\ell =1}^n \left( 
\int d \mes ( \gl_1 ) \dots d \mes ( \gl_{m_\ell} )
 \sum_{  \cG_\ell \in \pmb{\cC}_{m_\ell}} 
\prod_{ \{j,j'\} \in E(\cG_\ell)} \big(- \indc _{ \lambda_j \sim_\eps \lambda_{j'}}  \big) 
\right),
\end{align*}
where as previously  $\frac{k!}{m_1 ! \dots m_n ! }$ is the number of partitions  of~$\{1,\dots,k\}$  into~$n$ ordered components of cardinalities $m_1, \dots , m_n$ and~$\frac{1}{n!}$ removes the multiple counting due to the ordering. 
Using the definition \eqref{eq: phi} of $\gp_\eps$, we get 
\begin{align*}
\cZ^\eps_T (e^H)
= 1+  \sum_{k = 1}^\infty 
\sum_{n = 1}^\infty 1_{n \leq k}  \, \frac{1}{n!}
\sum_{m_1, \dots , m_n \geq 1 \atop   m_1 + \dots + m_n = k} 
\prod_{\ell =1}^n \left( \; \frac{1}{ m_\ell !}
\int d \mes ( \gl_1 ) \dots d \mes ( \gl_{m_\ell} )
 \, \varphi_\eps (\gl_1, \dots, \gl_{m_\ell}  )
\right).
\end{align*}
Choosing $T$ small enough, the sums are absolutely convergent thanks to \eqref{eq: bornes cluster expansion} so that they can be swapped  
\begin{align*}
\cZ^\eps _T (e^H)& = 1+   
\sum_{n = 1}^\infty   \frac{1}{n!}
\sum_{m_1, \dots , m_n \geq 1 } 
\prod_{\ell =1}^n \left( \; \frac{1}{ m_\ell !}
\int d \mes ( \gl_1 ) \dots d \mes ( \gl_{m_\ell} )
 \, \varphi _\eps(\gl_1, \dots, \gl_{m_\ell}  )
\right)\\
& = 1+   \sum_{n = 1}^\infty   \frac{1}{n!}
 \left(  \sum_{m \geq 1} \; \frac{1}{ m !}\int d \mes ( \gl_1 ) \dots d \mes ( \gl_m )
\, \varphi_\eps (\gl_1, \dots, \gl_m  )
\right)^n.
\end{align*}
This is the expansion of the exponential which can be inverted to recover 
\eqref{eq: cluster expansion}.

\bigskip

We turn now to the derivation of  the estimate \eqref{eq: bornes cluster expansion}
which relies on the specific structure of the microscopic dynamics   and more precisely on the geometry of the trajectories in $[0,T]$. 
We will use the geometric estimates  devised in \cite{BGSS2} (see also~\cite{BGSS3}).

\medskip

\noindent 
\emph{Case  $k =1$.}
We first consider   a single $\eps$-cluster path $\gl$ and prove the existence of a time $T_1>0$ and of a constant $c_0$ 
such that for~$T\leq T_1$
\begin{equation}
\label{eq: integration 1 foret}
\int d |\mes| (\gl) \; e^{\frac{10}{\beta}  |\gl| }  \leq  c_0 \mu_\eps\,,
\end{equation}
where  an additional term $e^{\frac{10}{\beta}  |\gl| }$ was added to inequality \eqref{eq: bornes cluster expansion} for later purposes.

\medskip
Using the definition \eqref{eq: distribution cluster path} of $\mes$ and summing over the size $n$ of the cluster path, one has
\begin{equation}
\label{eq: integration 2 foret}
\int d |\mes| (\gl) e^{\frac{10}{\beta}  |\gl| }   = 
\sum_{n \geq 1} \frac{ \mu_\eps^n }{n \, !} \int 
e^{ \frac{10}{\beta} n} \Big| \exp \big( H(\lambda) \big) \Big| \, \left(\prod_{i =1}^n f^0(z_i) \, 
  \right) \;  \indc_{\lambda  ~\hbox{\scriptsize $\eps$-cluster path of size} ~n  } \; d Z_n \,.
\end{equation}
Thanks to  assumption \eqref{eq: Gaussian bound} on the initial distribution $f^0$ and 
assumption \eqref{eq: hypotheses h} on $H$,
we have the following upper bound, for some constant~$C'_0 $ depending on~$C_0$ and~$\beta$
\begin{equation}
\label{eq: upper bound integrand}
 \Big| \exp \big( H(\lambda) \big) \Big|  \indc_{|\lambda|= n  } \, e^{ \frac{10}{\beta} n}  \, \prod_{i =1}^n f^0(z_i) 
\leq \Big(C'_0 e^{c_1 +\frac{10}\beta}\Big) ^n  \; \cM_{3\beta/4}^{\otimes n} (V_n) \, .
\end{equation}
Plugging (\ref{eq: upper bound integrand}) into (\ref{eq: integration 2 foret}), we see that  it is enough to show that there is $C>0$ such that  uniformly in $n \geq 2$
\begin{equation}
\label{eq: integration nu 2}
\int \cM_{\beta/4}^{\otimes n} (V_n) \;  \indc_{\lambda\hbox{ \scriptsize$\eps$-cluster path of size}~n  } \; d Z_n \;
\leq  n ! \frac{C^n  \, T^{n-1}}{ \mu_\eps^{n-1} } \;,
\end{equation}
and then to choose $T_1>0$ small enough to obtain  \eqref{eq: integration 1 foret} after summation over~$n$.
Note that only part of the Gaussian weight $\cM_{3\beta/4}^{\otimes n}$ in \eqref{eq: upper bound integrand} has been used in the above estimate, as we shall need an additional exponential decay later on  when dealing with the case~$k\geq 2$.

\begin{figure}[h] 
\centering
\includegraphics[width=3.2in]{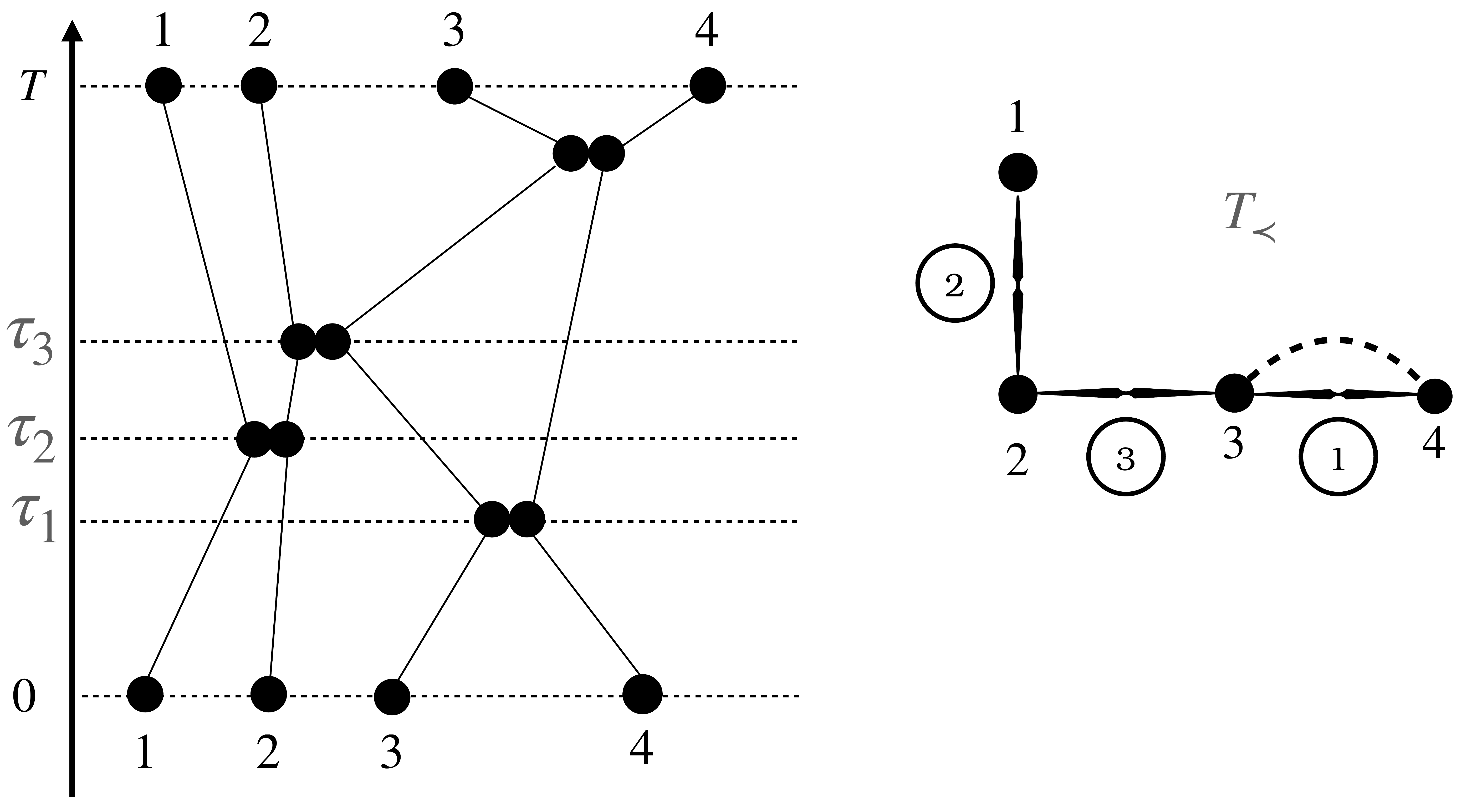} 
\caption{\small 
A single cluster path with 4 particles ($n = 4$) is depicted  on the left and the corresponding collision tree $\cT_\prec$ is represented on the right. The edges of the tree $\cT_\prec$ are ordered (circled numbers) according to the order of the collisions in the cluster path $\tau_1< \tau_2< \tau_3 \leq T$.
Recollisions may occur in the cluster path (as between particles $3$ and $4$ on the picture), but the corresponding (dashed) edge does not belong to $\cT_\prec$. The collisions associated with the graph edges are called clustering collisions.
}
\label{figure: 1 cluster path}
\end{figure}

\medskip
Let us now prove~(\ref{eq: integration nu 2}).
The constraint that $\lambda$ is an $\eps$-cluster path of size $n$ imposes that all the particles interact dynamically during the time interval $[0,T]$. We are going to record these collisions in an ordered tree~$\cT_\prec=  \{q_i, \bar q_i\}_{1\leq i \leq n -1}$ (see Figure \ref{figure: 1 cluster path}).
There can be more than $n-1$ collisions in the dynamics, but in order to retain a minimal structure and to end up with a tree~$\cT_\prec$, the collisions creating a cycle in the graph are not recorded. The collisions kept in the tree~$\cT_\prec$ will be called {\it clustering collisions}~:
 the first collision occurs between particles~$   q_1$ and~$\bar q_1$ at time~$\tau_1 \in [0,T]$, and the last collision is between~$   q_{n-1}$ and~$\bar q_{n-1}$ at time~$\tau_{n-1} \in (\tau_{n-2},T)$. 
In this way, an ordered graph recording the dynamical  interactions is built by following the flow of the hard sphere dynamics
in $[0,T]$. At intermediate times, the graph is made of several connected components,  and becomes completely connected at time~$\tau_{n-1}$.

The set of all ordered trees with $n$ vertices is denoted by $\pmb{\cT}^\prec_n$.
Thus summing over all the trees leads to
\begin{equation}
\label{eq: cluster path lambda1}
\indc_{\lambda ~ \hbox{\scriptsize $\eps$-cluster path of size} ~n  }= 
 \sum_{\cT_\prec \in \pmb{\cT}^\prec_n} 
 \indc_{\{ Z_n \in \cR_{\cT_\prec}  \}    } \, ,
\end{equation}
where $\cR_{\cT_\prec} $ is the set of initial configurations $Z_n$ with trajectories compatible with the ordered tree $\cT_\prec$.
By construction for fixed $\eps>0$ and for any given~$Z_n$  coding an $\eps$-cluster path, only one term is non zero in the  above sum.

\medskip

 For an admissible tree $\cT_\prec$, the relative positions, at the initial time, of the $i^\text{th}$  colliding 
 particles are  denoted by 
\begin{equation}
\label{eq: relative distance}
 \hat x_i:=x_{   q_i}-x_{\bar q_i}\, .
\end{equation}
Given the relative positions $(\hat x_j )_{j < i}$ and the velocities $V_n$, we consider a forward flow with clustering collisions at times $\tau_{ 1}< \dots < \tau_{i-1} < T$. By construction, $   q_i$ and $\bar q_i$ do not belong to the same connected component in the graph~$G_{i-1}$ defined as the graph with edges $\{   q_j,\bar q_j\}_{1 \leq j \leq i-1}$. 
We shall denote by $\cC_{   q_i}$ and by $\cC_{\bar q_i}$ the connected components associated 
with $   q_i,\bar q_i$ at time $\tau_{i-1}$.
Inside each connected component, the relative positions are fixed by the dynamical constraints, but the whole component can be translated so  that a free parameter remains. 
Therefore by varying $\hat x_i$ (moving rigidly the  connected components $\pmb{\cC}_{   q_i}, \pmb{\cC}_{\bar q_i}$), a  collision at time~$\tau_i \in ( \tau_{i-1}, T )$ 
 between $   q_i$ and $\bar q_i$ can be triggered. 
If the particles $  q_i$ and $\bar q_i$ move in straight lines, then the collision at time $\tau_i$ imposes a constraint at time $\tau_{i-1}$
\begin{equation}
\label{collision-eq}
\eps = \big| {\textbf x}^{\e}_{   q_i} ( \tau_i) - {\textbf x}^{\e}_{\bar q_i} ( \tau_i) \big|
= \big| {\textbf x}^{\e}_{   q_i} ( \tau_{i-1}) - {\textbf x}^{\e}_{\bar q_i} ( \tau_{i-1}) - (\tau_i - \tau_{i-1}) \big(  \textbf{v}^\eps_{  q_i} ( \tau_{i-1}^+) - \textbf{v}^\eps_{\bar q_i} ( \tau_{i-1}^+) \big) \big|\, .
\end{equation}
This says that the relative position at time  $\tau_{i-1}$ has to belong to a tube of direction 
$\textbf{v}^\eps_{   q_i} ( \tau_{i-1}^+) - \textbf{v}^\eps_{\bar q_i} ( \tau_{i-1}^+)$ with
diameter $\eps$ and length $| \textbf{v}^\eps_{  q_i}(\tau^+_{i-1}) -   \textbf{v}^\eps_{\bar q_i}(\tau^+_{i-1})| 	
\left( T - \tau_{i-1}\right)$ so that the collision occurs before time $T$.
 Since the trajectories inside each connected components  are fixed up to a global translation, this condition can be expressed in terms of the initial relative position~: $\hat x_i\in \cB_{\cT_\prec, i} (\hat x_{1}, \dots, \hat x_{i-1}, V_n )$.
 Thus, the measure of the set $\cB_{\cT_\prec, i}$ is bounded from above by (recall~$\mu_\eps= \eps^{-(d-1)}$)
\begin{equation}
\label{eq: borne B T prec i}
|\cB_{\cT_\prec, i} | 
\leq \frac{C}{\mu_\eps} | \textbf{v}^\eps_{  q_i}(\tau^+_{i-1}) -   \textbf{v}^\eps_{\bar q_i}(\tau^+_{i-1})| 	
\; \int^T_{\tau_{i-1}} d \tau_i \; .
\end{equation}
If  the particle $  q_i$ (resp. $\bar q_i$) has been deflected during $[\tau_{i-1},\tau_i[$  (by recollisions with particles in the connected component of $\pmb{\cC}_{  q_i}$ (resp. $\pmb{\cC}_{\bar q_i}$)) then  one has to decompose the trajectories into a union of tubes (as in Chapter 8 of \cite{BGSS2}) and  an estimate as~(\ref{eq: borne B T prec i})  can be recovered.
Summing over all the possible pairs of particles and using the fact that collisions preserve the kinetic energy,  we get after a Cauchy-Schwarz inequality
\begin{equation}
\label{eq: borne somme qi}
\sum_{  q_i,\bar q_i}|\cB_{\cT_\prec, i} |  \leq \frac{C}{\mu_\eps} \left( {\mathcal E}_\lambda +  n \right) n 	
\; \int^T_{\tau_{i-1}} d \tau_i \; ,
\end{equation}
where ${\mathcal E}_\lambda:=V_n^2/2 $ is the total kinetic energy of the particles in the cluster path.  
Iterating the previous estimates, we get from Fubini's theorem
\begin{align}
\sum_{\cT_\prec \in \pmb{\cT}^\prec_n }
\int d\hat X_{n-1} 
\prod_{i=1}^{n-1}  \indc_{\cB_{\cT_\prec, i }} & \leq 
\sum_{\cT_\prec \in \pmb{\cT}^\prec_n } \int d\hat x_1
 \indc_{\cB_{\cT_\prec, 1 }} \int d\hat x_{ 2}\,  \dots\int d\hat x_{n-1} \indc_{\cB_{\cT_\prec, n-1}}
 \nonumber \\
 &  \leq \left( \frac{C}{\mu_\eps}\right)^{n-1}  \left({\mathcal E}_\lambda + n \right)^{n-1} n^{n-1}
 \int_0^T d\tau_1\, \dots \int_{\tau_{n-2} }^T d \tau_{n-1}
 \leq 
 \left( \frac{C}{\mu_\eps}\right)^{n-1}  \left({\mathcal E}_\lambda + n \right)^{n-1} n^{n-1} \frac{T^{n-1}}{(n-1)!}  ,
 \label{eq: Bi taille}
\end{align}
where the last inequality follows by integrating the ordered times.
Furthermore,  for any $K,N$   
\begin{equation}
\label{eq: inegalite exponentielle}
\sup_{ V \in \bbR^{d N}} \Big\{ \exp \big( - \frac{\beta}{4} |V|^2 \big)  \; ( |V|^2 + K)^N  \Big\} \leq C_\beta^N e^{\frac \beta 4K} \; N^N.
\end{equation}
Using the Gaussian integration in \eqref{eq: integration nu 2}, we deduce from the previous inequality that the term $\left({\mathcal E}_\lambda+ n \right)^n$ leads to another factor of order $n^n$ which 
is (up to a factor $C^n$) of the same order as $n !$. 
 Recalling~(\ref{eq: cluster path lambda1}), this completes \eqref{eq: integration nu 2} and thus \eqref{eq: integration 1 foret}.

 \begin{Rmk} \label{rk : jacobian} Let us call the  \emph{root} of~$\lambda$ the position~$y:=\displaystyle \frac1n \sum_{i=1}^nx_i$ of the center of mass of   the particles in~$\lambda$ at time~0,   where $n = |\gl|$.
The change of variables $X_n \mapsto ( \hat X_{n-1}, y)$ has unit Jacobian. Moreover,  since the $\eps$-cluster path $\lambda$ is invariant by translations, its root  is not constrained by the collision conditions.
 \end{Rmk}

\bigskip

\noindent 
\emph{Case $k >1$.}
We derive now \eqref{eq: bornes cluster expansion} for a non trivial $\eps$-aggregate with $k$ $\eps$-cluster paths: let us prove that for~$T$ small enough (smaller than $T_1$ (see \eqref{eq: integration 1 foret}) and independent of~$k$)
\begin{equation*}
\int d | \mes| ( \gl_1 ) \dots d |\mes| (  \gl_k ) \; \varphi_\eps (\gl_1, \dots, \gl_k )  
\leq   \mu_\eps \;  k!   \, C^k \,  \big(  T + \eps \big)^{k-1}.
\end{equation*}
For this, we are going to use the inequality \eqref{eq: integration 1 foret} which evaluates the constraints, 
in each cluster path $\gl_\ell$, on the coordinates of the particles $Z^{(\ell)} _{|\gl_\ell|}$ at time 0.
Further dynamical constraints are added by the function $\varphi_\eps (\gl_1, \dots, \gl_k )$ defined in \eqref{eq: phi}. Viewing the overlaps between the $\eps$-cluster paths as the edges of a graph  with $k$ vertices indexing the $\eps$-cluster paths $\gl_1, \dots, \gl_k$, then the alternating sums defining $\varphi_\eps (\gl_1, \dots, \gl_k)$ can be bounded from above by 
the so called {\it tree inequality} 
\begin{equation}
\label{eq: overlap tree}
\big| \varphi_\eps (\gl_1, \dots, \gl_k) \big| = \Big| \sum_{\cG \in \pmb{\cC}_k} \prod_{\{ j,j'\}  \in E (\cG)} (-\indc_{\lambda_ j \sim_\eps \lambda_{j'}}) \Big| 
\leq  
\sum_{\cT^{ov}  \in  \pmb{\cT}_k}  \  \  \prod_{\{ j,j'\}  \in E (\cT^{ov} )} \indc_{\lambda_ j \sim_\eps \lambda_{j'}}\,,  
\end{equation}
where the sum is restricted to  (non-ordered) trees. 
The tree inequality \eqref{eq: overlap tree} which allows to control $\varphi_\eps$  follows from a standard combinatorial argument  (see \cite{penrose1967convergence}) which will not be recalled here.
Note that in \eqref{eq: overlap tree}, several graphs can be compatible with the same family of $\eps$-cluster paths  $ (\gl_1, \dots, \gl_k)$
  so that several trees $\cT^{ov}$ may contribute to the sum.

\medskip

Since the particle trajectories are unchanged by the overlaps, it is not needed to proceed as for the collisions within the $\eps$-cluster paths and to prescribe an order related to the dynamical overlaps on the edges of $\cT^{ov}$.
Thus, for each $\cT^{ov}$, we  have more flexibility in choosing the integration variables and we can order the edges following for instance the depth of vertices in the graph (see \cite{BGSS2} for details). We then examine successively the $k-1$ overlap constraints imposed on the $\eps$-cluster paths.
Denote by $  \gl_{  p_i} ,   \gl_{\bar p_i}$ the  $\eps$-cluster paths involved in the $i^{\text {th}}$ overlap and by $    q_i $ 
and $\bar  q_{i} $ the two overlapping particles.
The constraint imposed by the $i^{\text {th}}$ overlap leads to 
a condition on  the  relative position of the two cluster paths   at the initial time  (see \eqref{eq: relative distance})
\begin{equation*}
\hat y_i:=  y_{  p_i} -y_{  \bar p_i }\,,
\end{equation*}
recalling the notation for the root introduced in Remark~\ref{rk : jacobian}.
Indeed, fixing the velocities and the relative positions in each cluster path, one has to evaluate the condition for the set 
$$\big\{ \hat y_i + (x_{  q_i}(t) -  y_{  p_i}) - (x_{\bar q_i}(t) -  y_{\bar p_i})  ; \quad t \leq T \big\}$$ 
to intersect a ball of radius $\eps$ around the origin. 
Note that, contrary to the collisions, the overlap may occur at the initial time or dynamically. 
This condition on $\hat y_i$ is coded by the set $\cB_{\cT^{ov}, i}$  with measure bounded by
\begin{equation}
\label{eq: cout overlap}
|\cB_{\cT^{ov}, i}| 
\leq C \eps^d 
+ \frac{C}{\mu_\eps} \int_0^T ds \; | \textbf{v}^\eps_{  q_i}(s) - \textbf{v}^\eps_{\bar q_i}(s)| \, .
\end{equation}
We stress the fact that $\eps^d$ corresponds to the cost of an overlap at time 0 which is much smaller than the cost $\frac{1}{\mu_\eps} = \eps^{d-1}$ of a dynamical overlap.
This fact will be used in the Boltzmann-Grad limit to neglect the overlaps occurring at the initial time.

Denoting as previously by $|\gl |$ the cardinality of an $\eps$-cluster path and summing over all the possible particles in the $\eps$-cluster paths
$  \gl_{  p_i} ,   \gl_{\bar p_i}$, we get by a Cauchy-Schwarz inequality
\begin{align}
\sum_{  q_i, \bar q_i } |\cB_{\cT^{ov}, i}|  
& \leq C\eps^d \; |\gl_{\bar p_i}| \; |  \gl_{ p_i} | 
+ \frac{C}{\mu_\eps} \int_0^T ds \; 
\left(  |\gl_{   p_i}| \; \sqrt{ |   \gl_{\bar p_i} | }  \sqrt{ \cE_{\gl_{\bar p_i}}}
+ 
|\gl_{\bar p_i}| \; \sqrt{ |\gl_{   p_i}|} \sqrt{ \cE_{\gl_{  p_i}} }\right)
\nonumber \\
& \leq C\eps^d \; | \gl_{\bar p_i}| \; |\gl_{  p_i}| 
+ \frac{C}{\mu_\eps} T  \; 
 \Big( \frac{\beta}{4}  \cE_{\gl_{\bar p_i}} + \frac{4}{\beta} |\gl_{\bar p_i}| \Big)  
 \Big(  \frac{\beta}{4} \cE_{\gl_{  p_i}} + \frac{4}{\beta} |\gl_{ p_i}| \Big) 
\leq  \frac{C}{\mu_\eps} (T  + \eps )\; 
 \Big( \frac{\beta}{4}  \cE_{\gl_{\bar p_i}} + \frac{4}{\beta} |\gl_{\bar p_i}| \Big) 
\Big( \frac{\beta}{4}  \cE_{\gl_{  p_i}} + \frac{4}{\beta} |\gl_{ p_i}| \Big)   ,
\label{eq: overlap 2 forets}
\end{align}
where we used in the second inequality that the total kinetic energy  $\cE_{\gl}$ of the particles in an $\eps$-cluster path $\gl$ 
is constant in time.
Thus \eqref{eq: overlap 2 forets} measures the cost of the overlap coded by the edge $\{\gl_{\bar p_i}, \gl_{  p_i}\}$ in the tree $\cT^{ov}$.

Let $\cB_{\cT^{ov}}$ be the set representing all the conditions imposed by the overlaps.
Given $\{ \cE_{\gl_1}, \dots, \cE_{\gl_k}\}$ the energies of all the $\eps$-cluster paths, the measure of 
$\cB_{\cT^{ov}}$   with respect to the relative positions~$(\hat y_i)_{i \leq k-1}$ of the $\eps$-cluster paths 
  is obtained by multiplying the contributions 
  \eqref{eq: overlap 2 forets} for each edge of the tree $\cT^{ov}$ 
\begin{equation*}
|\cB_{\cT^{ov}}| \leq \left( \frac{C}{\mu_\eps}\right)^{k-1} \;  (T  + \eps )^{k-1} \; \prod_{i = 1}^k 
 \Big( \frac{\beta}{4}  \cE_{\gl_i} + \frac{4}{\beta} |\gl_i| \Big)^{d_i} ,
\end{equation*}
where $d_i$ stands for the degree of the vertex $\gl_i$ in the tree $\cT^{ov}$.
There are $( k - 2 )! / \prod_i  ( d_i-1)!$
 trees of size $k$ with specified vertex degrees (see e.g.\,Lemma 2.4.1 in \cite{BGSS2}).
Thus summing over all the trees $\cT^{ov}$, we get 
\begin{align*}
\sum_{\cT^{ov} \in  \pmb{\cT}_k} |\cB_{\cT^{ov}}| 
& \leq \left( \frac{C}{\mu_\eps}\right)^{k-1} \; (T  + \eps )^{k-1} ( k - 2 )! 
\sum_{d_1, \dots d_k \atop d_1 + \dots d_k = 2 k -2} \; 
\prod_{i =1}^k  \frac {\Big( \frac{\beta}{4}  \cE_{\gl_i} + \frac{4}{\beta} |\gl_i| \Big)^{d_i} }{ \left( d_i -1\right)! }\\
& \leq \left( \frac{C}{\mu_\eps}\right)^{k-1} \; (T  + \eps )^{k-1} ( k - 2 )! \;
\prod_{i =1}^k \left(  \frac{\beta}{4}   \cE_{\gl_i} + \frac{4}{\beta}  |\gl_i| \right) 
\; \exp \left(  \frac{\beta}{4}   \cE_{\gl_i} + \frac{4}{\beta}  |\gl_i| \right),
\end{align*}
where the constraint on the degrees is released in the last inequality to recover the exponential.

Using the Gaussian weights $\cM_{\beta/2}^{\otimes  |\gl_i| }$ from the initial measure as well as  the inequality  \eqref{eq: inegalite exponentielle}, we obtain an upper bound for the overlaps of the form 
\begin{align}
\sum_{\cT^{ov} \in  \pmb{\cT}_k} |\cB_{\cT^{ov}}|   \; \prod_{i =1}^k \cM_{\beta/2}^{\otimes  |\gl_i| } (Z_{ |\gl_i|})
\leq 
 \frac{C^k}{\mu_\eps^{k-1}}  \, k ! \, (T  + \eps )^{k-1} 
\prod_{i =1}^k  e^{\frac{10}{\beta}  |\gl_i| }.
\label{eq: final overlap}
\end{align}
Once the dynamical constraint on the overlaps has been taken into account, the 
contributions of the cluster paths are independent and can be estimated by \eqref{eq: integration 1 foret}
for $T$ small enough (independently of~$k$)
\begin{equation}
\label{eq: integration nu 3}
\prod_{i =1}^k  \int  d | \mes| ( \gl_i ) \, e^{\frac{10}{\beta}  |\gl_i| }\; \leq  c_0^k  \mu_\eps^k .
\end{equation}
Combining \eqref{eq: overlap tree}, \eqref{eq: final overlap} and \eqref{eq: integration nu 3}, we deduce that \begin{equation}
\int d | \mes| ( \gl_1 ) \; \dots  \;  d |\mes| ( \gl_k ) \; \varphi _\eps (\gl_1, \dots, \gl_k )  
\leq   \mu_\eps \;  k!  \; C^k \, (T  + \eps )^{k-1} .
\end{equation}
This completes  the proof of \eqref{eq: bornes cluster expansion} for a value of $T>0$ small enough so that 
\eqref{eq: integration 1 foret} holds. Proposition \ref{prop: cluster expansion} is proved.

%
%
%

\subsection{Measure concentration properties}

We are going to deduce now some consequences of the cluster expansion 
for fixed (but large)  $\mu_\eps$.
 The exponential moments~$\gL^\eps_T$ introduced in  \eqref{eq: log modified partition function}
 encode the statistics of the cluster paths.
 Choosing $T>0$ as in Proposition \ref{prop: cluster expansion}, 
 we know that the functional~$\gL^\eps_T (e^{uH})$ is analytic with respect to~$u$
 and its derivatives at $u =0$ are uniformly controlled in $\mu_\eps$.
 This is an important property as the derivatives of the functional are related to 
physical quantities.
Indeed considering for instance the first derivative at 0 of~$u \in \bbR \mapsto \gL^\eps_T (e^{u \, H})$, we recover the expectation of the empirical measure (using the notation of \eqref{eq: modified empirical general})
\begin{equation}
\label{eq: average Pi}
\bbE_\eps \left[  \Pi^\eps_{[0,T]}(H)   \right] 
= \partial_u \, \gL^\eps_T (e^{u \, H}) \Big|_{u =0}
= \frac{1}{\mu_\eps} 
\sum_{k \geq 1}\frac{1}{k!}  
\int d \mesflat ( \gl_1 ) \dots d \mesflat ( \gl_k ) \; \varphi_\eps (\gl_1, \dots, \gl_k ) 
\left( \sum_{\ell =1}^k    H  \big( \gl_\ell \big) \right).
\end{equation}
Using the symmetry between the clusters, 
we recover the distribution of a typical $\eps$-cluster path
\begin{equation}
\label{eq: esperance tanakae}
\bbE_\eps \left[  \Pi^\eps_{[0,T]}(H)   \right]  = \int d\tanakae (\gl_1) H(\gl_1),
\end{equation}
where \begin{equation}
\label{eq: def tanakae}
d\tanakae (\gl_1) = \frac{1}{\mu_\eps} 
\left( d\mesflat ( \gl_1 )  + \sum_{k \geq 2}\frac{1}{(k-1)!}  d\mesflat ( \gl_1 ) 
\int   d \mesflat ( \gl_2 )  \dots d \mesflat ( \gl_k  ) \; \varphi_\eps (\gl_1, \dots, \gl_k )  \right).
\end{equation}

Taking twice the derivative leads to the variance 
\begin{equation}
\bbE_\eps \left[  \left(\Pi^\eps_{[0,T]}(H)    - \bbE_\eps \left[  \Pi^\eps_{[0,T]}(H)   \right] \right)^2 \right]
= \bbE_\eps \left[  \Pi^\eps_{[0,T]}(H)^2   \right] - \bbE_\eps \left[  \Pi^\eps_{[0,T]}(H)   \right]^2
= \frac{1}{\mu_\eps} \partial^2_u \, \gL^\eps_T (e^{u \, H}) \Big|_{u =0} .
\end{equation}
As a corollary of Proposition \ref{prop: cluster expansion},
the second derivative is uniformly bounded with respect to~$\mu_\eps$ \begin{align}
\partial^2_u \, \gL^\eps_T (e^{u \, H}) \Big|_{u =0}
& = \frac{1}{\mu_\eps} \sum_{k \geq 1}\frac{1}{k!}  
\int d \mesflat ( \gl_1 ) \dots d \mesflat ( \gl_k ) \; \varphi_\eps (\gl_1, \dots, \gl_k )   
\left( \sum_{\ell =1}^k  H  \big( \gl_\ell \big) \right)^2\\
& \leq   \frac{1 }{\mu_\eps} \sum_{k \geq 1}\frac{1}{k!}   \int d |\mesflat| ( \gl_1 ) \dots d |\mesflat |( \gl_k ) \; \varphi_\eps (\gl_1, \dots, \gl_k ) \left(\sum_{\ell =1}^k c_1 | \gl_\ell |+c_2\cE_{\gl_\ell} \right) ^2 = O(1)  \, , \nonumber
\end{align}
for test functions $H$ satisfying~(\ref{eq: hypotheses h}).

This implies that the covariance vanishes in the Boltzmann-Grad limit so that the empirical measure concentrates to its mean
\begin{equation}
\label{eq: covariance Pi}
\bbE_\eps \left[ \left(  \Pi^\eps_{[0,T]}(H)   - \bbE_\eps \left[  \Pi^\eps_{[0,T]}(H) \right]  \right)^2 \right]
= O \left( \frac{ 1 }{\mu_\eps}  \right) .
\end{equation}

\bigskip

Applying this result to functions $H(\lambda) =\displaystyle  \sum_{i =1}^{|\lambda|} h   \big( {\textbf z}_i ([0,T]) \big)$
as in \eqref{eq: hypotheses h bis}, we deduce in particular the concentration of the empirical measure \eqref{eq: empirical general} to its mean
\begin{equation}
\label{eq: covariance pi}
\bbE_\eps \left[ \left(  \pi^\eps_{[0,T]}(h)   - \bbE_\eps \left[  \pi^\eps_{[0,T]}(h) \right]  \right)^2 \right]
= O \left( \frac{\| h \|_\infty^2}{\mu_\eps}  \right) .
\end{equation}
By taking further derivatives, one can recover all the cumulants studied in~\cite{BGSS2}  and show that the $L^1$-norm of the (unrescaled) cumulant of order $n$ decays as $O \left( \mu_\eps^{1-n}  \right)$.
This result was already obtained in \cite{BGSS2} (see Theorem 4 therein). Notice however that the series expansion for $\widetilde \gL^\eps_T (e^h)$ is derived in \cite{BGSS2} by applying a cluster expansion on  the Duhamel representation of the correlation functions and the terms of the series are described by pseudo-trajectories instead of physical trajectories (see Eq.s (4.4.7) and (4.4.1) in  \cite{BGSS2}).


\section{The Boltzmann-Grad limit}
\label{sec: Boltzmann-Grad limit}

The series expansion \eqref{eq: cluster expansion} in Proposition \ref{prop: cluster expansion}
provides a representation of the partition function $\log  \cZ^\eps_T (e^H)$ in terms of generalized {\it dynamical interactions} (collisions and overlaps) of microscopic trajectories.  The estimates derived in the previous section hold uniformly with respect to $\mu_\eps$ (large enough).
In this section, we are going focus on the kinetic limit~$\mu_\eps\to \infty$.
We are going to show that the structure of these interactions simplifies in that  limit, 
providing a simpler (but singular) expansion.

\subsection{Discarding recollisions and non  minimal overlaps}

\begin{figure}[h] 
\centering
\includegraphics[width=3.5in]{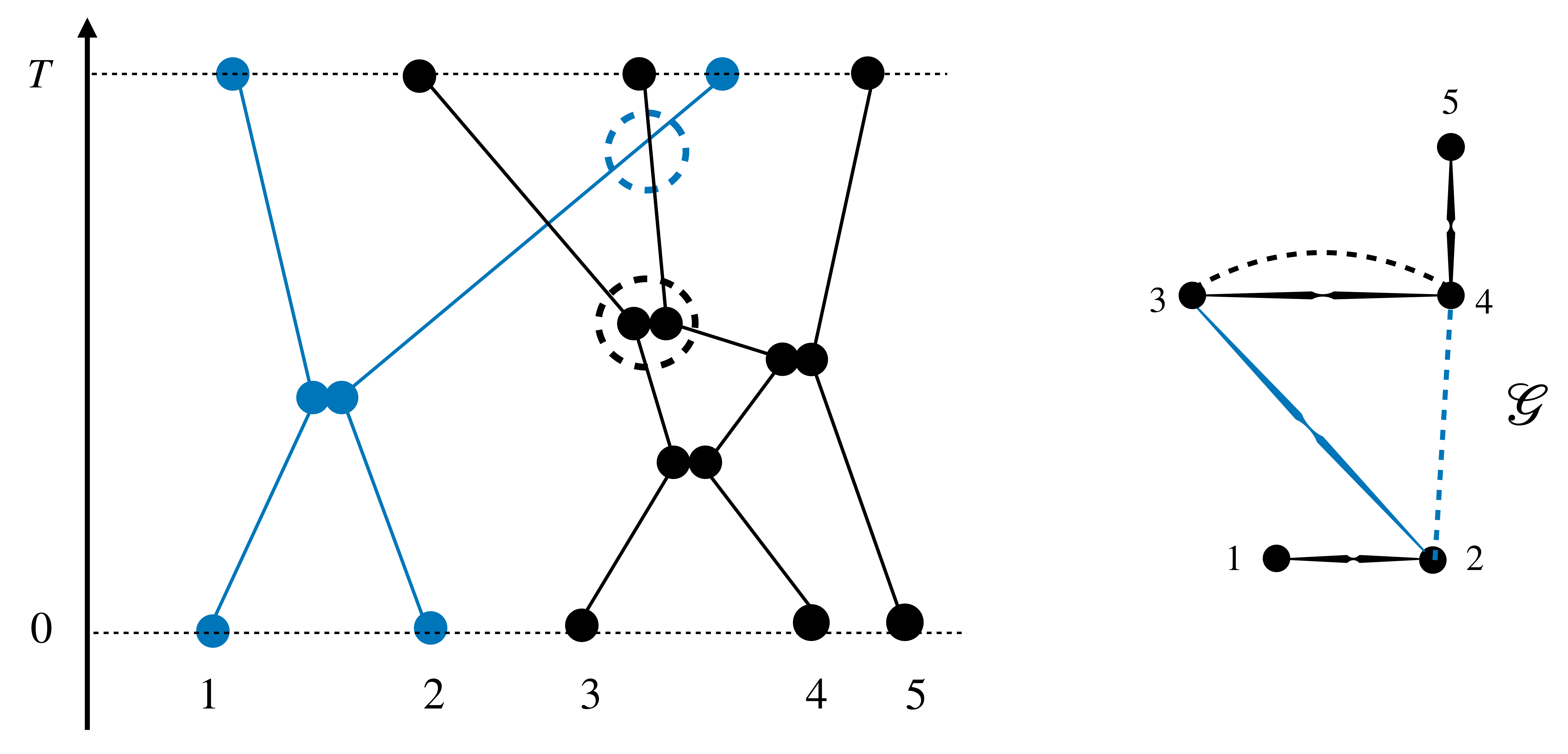} 
\caption{\small On the left, two overlapping $\eps$-cluster paths are represented with different colors.  The dynamical interaction graph is depicted on the right with black edges for  the collisions and blue edges for the overlaps. The minimal dynamical graph is shown in solid lines, instead the recollisions and multiple overlaps are represented by dashed edges.}
\label{figure: recollisionreoverlap}
\end{figure}

We first consider the dynamical interactions  (collisions) within an $\eps$-cluster path $\gl$  and show that,   in the Boltzmann-Grad limit, the only relevant trajectories have exactly $|\gl |-1$ collisions.
Recall that in  this $\eps$-cluster path $\lambda$, all the particles interact dynamically  in the sense of Definition \ref{def: cluster paths}. 
All these interactions can be recorded in a  {\it dynamical interaction graph} $\cG$ with $ |\lambda| $ vertices labelled by the particles and edges corresponding to a collision between two particles. By definition of the $\eps$-cluster path $\gl$, the graph $\cG$ is connected and it may have cycles or multiple edges (see  the black cluster path in Figure \ref{figure: recollisionreoverlap}). 

These cycles correspond to \textit{recollisions} in the hard sphere dynamics starting from the initial configuration $Z_{|\lambda | } $. 
We stress the fact that this notion of recollision slightly differs from the interpretation of recollisions used in the Duhamel representation \cite{Cercignani_Illner_Pulvirenti, Spohn_fluctuations, Pulvirenti_Simonella, BGSS2}.

In the Boltzmann-Grad limit,  the only relevant graphs will be trees, i.e.\,minimally connected graphs.
To discuss the convergence in this limit, it is then useful to introduce a truncated measure
which (with respect to \eqref{eq: distribution cluster path}) forbids the recollisions. We say that 
an $\eps$-cluster path $\lambda$ is minimal, or (for brevity) a {\it  min $\eps$-cluster path}, if its dynamical interaction graph is a tree. 
We define the following truncated integration measure~:
\begin{equation}
\label{eq: reduced distribution forest}
d\mesh0  ( \gl ) := \frac{\mu_\eps^{| \gl |}}{| \gl | \, !}\,  \indc_{ {\lambda} \text{ \scriptsize min $\eps$-cluster path on $[0,T]$ } }\, F^{(H)} (\gl)
\;  d Z_{|\lambda|}\;.
\end{equation}

In the proof of Proposition \ref{prop: cluster expansion}, 
Estimate \eqref{eq: integration 1 foret}, recalled below, 
\begin{equation*}
\int d |\mes| (\gl) \, e^{\frac{10}{\beta}  |\gl| }   \leq  C_0 \mu_\eps
\end{equation*}
was derived by showing that the collisions between particles in an $\eps$-cluster path can be indexed by a tree.
This tree records the minimal amount of dynamical constraints and the recollisions add more constraints which can be controlled  by the geometric  estimates derived in \cite{BGSS3}
(see Eq (5.12) and Appendix B), leading to 
\begin{equation}
\label{eq: correction lambda 1}
\Big| \int d\mesh0 ( \gl ) \, e^{\frac{10}{\beta}  |\gl| }  
 - \int d \mes ( \gl ) \, e^{\frac{10}{\beta}  |\gl| }   \Big|
\leq C \,  \mu_\eps \;  \eps ^\alpha
\end{equation}
for   $\alpha <1$, in dimension~$d\geq 3$ and for~$T$ small enough as in Proposition \ref{prop: cluster expansion}.  In dimension 2, similar geometric arguments provide the same bound.

\bigskip

The second  type of dynamical interactions are the overlaps, introduced in Definition \ref{def:overlap}.
We consider now $k$ overlapping cluster paths $\gl_1, \dots, \gl_k$. 
The combinatorial factor $\varphi_\eps (\gl_1, \dots, \gl_k )$ has been estimated in \eqref{eq: overlap tree}, recalled below, as an upper bound on trees with edges between two overlapping $\eps$-cluster paths
\begin{equation*}
\big| \varphi _\eps(\gl_1, \dots, \gl_k) \big| 
\leq  
\sum_{\cT^{ov}  \in  \pmb{\cT}_k}  \prod_{\{ j,j'\}  \in E (\cT^{ov} )} \indc_{\lambda_ j \sim_\eps  \lambda_{j'}}.
\end{equation*}
 If the graph recording all the overlaps of $(\gl_1, \dots, \gl_k)$ has cycles, then several trees $\cT^{ov}$ will contribute to the sum above. 
But we can actually show that cycles (see Figure \ref{figure: recollisionreoverlap}) can be neglected in the limit and that typically the constraint $\varphi _\eps (\gl_1, \dots, \gl_k)$ is compatible with a single tree. 
Furthermore, as noted in the comment after \eqref{eq: cout overlap}, the overlaps occurring initially have a much smaller cost than the dynamical overlaps. Thus they can also be neglected in the Boltzmann-Grad limit.

This means that  also in the case of overlaps, the only relevant dynamical interaction graphs are trees. When $\{\gl_1, \dots, \gl_k\}$ has such a minimal interaction graph, we say that it is {\it a minimal $\eps$-aggregate}.
By definition, the aggregate function $\varphi_\eps$ restricted to minimal $\eps$-aggregates  takes  the value $(-1)^{k-1}$ and is defined by 
\begin{equation}
\label{eq: overlap no cycle}
\varphi_\eps^{\hbox{\scriptsize min}}  (\gl_1, \dots, \gl_k) 
:=  (-1)^{k-1}   1_{\{\gl_1, \dots, \gl_k \} \hbox{\scriptsize  min $\eps$-aggregate}}\,  .
\end{equation}

Combining the  proof of Proposition \ref{prop: cluster expansion} and 
the geometric estimates derived in \cite{BGSS3}, one can show that 
the overlaps forming cycles do not contribute in the Boltzmann-Grad limit.
As the recollisions can also be neglected thanks to \eqref{eq: correction lambda 1},
we finally obtain that the minimally connected graphs provide the leading contribution 
to the cluster expansion series.
\begin{Prop}
\label{Prop: no cycle terms}
Let $T>0$ be the convergence time obtained in Proposition \ref{prop: cluster expansion}.
There   are constants $C,      \alpha \in (0,1)$  such that uniformly in $\eps$ small enough
and for $k\geq 1$
\begin{equation}
\label{eq: simplified terms cluster expansion}
\left|   
\int d\mesh0 ( \gl_1 ) \dots d\mesh0 ( \gl_k ) \; \varphi_\eps^{\hbox{ \rm\scriptsize min} }(\gl_1, \dots, \gl_k )   
-
\int d \mes ( \gl_1 ) \dots d \mes ( \gl_k ) \; \varphi_\eps (\gl_1, \dots, \gl_k )   
\right| 
\leq   \mu_\eps \;   \eps^{ \alpha}  \; k! \;  C^k  (T+\eps)^{k-1}.
\end{equation}
As a consequence,
\begin{equation}
\label{eq: simplified expansion}
\left| 
\log \cZ_T^\eps (e^H)  - \sum_{k \geq 1}\frac{1}{k!}  
\int d\mesh0 ( \gl_1 ) \dots d\mesh0 ( \gl_k ) \; \varphi_\eps^{\hbox{ \rm\scriptsize min} } (\gl_1, \dots, \gl_k )
\right|
\leq  C \,  \mu_\eps \;   \eps^{ \alpha}   .
\end{equation}
\end{Prop}

\begin{Rmk}\label{global parametrization - rmk} The term  $\displaystyle \int d\mesh0 ( \gl_1 ) \dots d\mesh0 ( \gl_k ) \; \varphi_\eps^{\hbox{ \rm\scriptsize min} } (\gl_1, \dots, \gl_k )$ in \eqref{eq: simplified terms cluster expansion} 
corresponds to   minimally connected dynamical graphs, and
can be rewritten differently by treating collisions and overlaps in a more symmetric way.

Given $\{ \gl_1, \dots, \gl_k \}$ a set of overlapping $\eps$-cluster paths, the corresponding particle configuration will be linked by $| \gl_\ell | -1$ collisions in each cluster path $\gl_\ell$ and $k-1$ overlaps between cluster paths.
In particular, if $n$ stands for the total number of particles in $\{ \gl_1, \dots, \gl_k \}$, 
we say that there are 
$
n-1 = \sum_{\ell =1}^k (| \gl_\ell | -1) + k-1
$
``clustering conditions'' (generalizing the clustering collisions introduced at the level of \eqref{eq: cluster path lambda1}).

Ordering all these clustering conditions according to the forward flow, we index them by a single signed ordered tree $\overline {\cT_\prec}$ such that the $n$ particles form the vertices and each edge $e$ has a sign $s_e = +1$ if it records a collision or $s_e = -1$ if it records an overlap. 
Thus instead of decomposing the particles into overlapping $\eps$-cluster paths, 
one can choose globally a set of $n$ particles   whose trajectories are denoted by ${\textbf Z}^\eps_n :={\textbf Z}^\eps_n ([0,T])$ and signed ordered trees~$\overline {\cT_\prec}$ coding the clustering conditions. 
By  Fubini's theorem, we deduce that 
\begin{equation*}
\sum_{k \geq 1}\frac{1}{k!}   \int d\mesh0 ( \gl_1 ) \dots d\mesh0 ( \gl_k ) \; \varphi_\eps^{\hbox{ \rm\scriptsize min} } (\gl_1, \dots, \gl_k )   
 = \sum_{n =1}^\infty \frac{\mu_\eps^n }{ n !} \sum_{\overline {\cT_\prec} \in \pmb{\cT}_n^{\prec, \pm}}
 \;  \Big( \prod_{e \in \overline {\cT_\prec} } s_e \Big) \; \int d Z_n \; F^{(H)} ({\textbf Z}^\eps_n) 
\; \indc_{  {\textbf Z}^\eps_n \   \hbox{\rm \scriptsize compatible with } \overline {\cT_\prec} }  
\end{equation*}
where $\pmb{\cT}_n^{\prec, \pm}$ is the set of signed ordered   trees with $n$ vertices
indexing the particle trajectories with initial data $Z_n$. 
\end{Rmk}

\subsection{Asymptotics of the partition function}
\label{sec: Asymptotics free energy}

We are going to use Proposition \ref{Prop: no cycle terms} to compute the asymptotics of the partition function when $\mu_\eps$ tends to infinity.
In the Boltzmann-Grad limit, the particle trajectories in the series expansion \eqref{eq: cluster expansion} become singular, but a limiting structure can be identified using a suitable parametrization of the $\eps$-cluster paths.

\bigskip

We first consider the dynamical interactions within a min $\eps$-cluster path $\gl$ of size $n$, and assume that the collisions are prescribed by the
ordered tree $ \cT_\prec$.  For fixed~$\eps>0$, the collision condition associated with the edge $e =\{i, j\}$ takes the form 
\begin{equation}
\label{eq: omega collision}
 \omega_{e } := \frac{ {\textbf  x}_i ^\e(\tau_e) - {\textbf  x}_j^\e (\tau_e  )} {\eps}  \in {\mathbb S}^{d-1} \,,
\end{equation}
denoting by $\tau_e$ the collision  time.
Recall that before the collision, the particles~$i,j$ are connected to two distinct components of the dynamical graph which can move rigidly   with  the positions of $i$ and~$j$ at time zero.
We know that the coordinates $\big({\textbf  x}_i ^\e (\tau_e^-) - x_i, {\textbf  v}_i ^\e (\tau_e^-) \big)$  and $\big( {\textbf  x}_j ^\e (\tau_e^-) - x_j, {\textbf  v}_j ^\e (\tau_e^-) \big)$  are fixed by the previous collisions (inside each dynamical component associated with $i$ and $j$). Then using (\ref{collision-eq}), we define the local change of variables
\begin{equation}
\label{eq: change variables}
\hat x_e := x_j - x_i \in\bbT^d \mapsto(\tau_{e } ,  \omega_{e }) \in  [0,T] \times  {\mathbb S}^{d-1}
\end{equation}
with  inverse Jacobian determinant $ \mu_\eps^{-1}\big( ( {\textbf  v}_i ^\e(\tau_e^-) - {\textbf  v}_j ^\e(\tau_e^-) )\cdot \omega_{e}\big)_+$. 
  This provides  the   identification of measures
\begin{equation}
\label{eq: identification measures}
 \mu_\eps \, dy \,dv_i\, d\hat x_e\,dv_j= dy  \, dv_i\,dv_j\,d\tau_{e} \; d\omega_{e }
\big(   ( {\textbf  v}_i ^\e(\tau_e^-) - {\textbf  v}_j ^\e (\tau_e^-) )\cdot \omega_e \big)_+ \,.
\end{equation}
Applying iteratively this change of variables for each edge of the  collision tree  $\cT_\prec$, the  $\eps$-cluster path~$\gl$    can be recovered by the following parameters~:
\begin{itemize}
\item $y$ the root of the $\eps$-cluster path (the center of mass of the~$n = |\gl|$ positions at time~0)
\item $(v_1, \dots, v_n)$ the particle velocities at the initial time,
\item $\Omega_{n -1} := \left(\omega_{e }\right)_{e \in E(\cT_\prec)}$ representing the collision angles,
\item $\Theta_{n -1} := \left(\tau_{e }\right)_{e \in E(\cT_\prec)}$ representing the collision times.
\end{itemize}
We  define the  singular measure 
\begin{align}
\label{eq: sing measure}
d\mu^\eps_{{\rm sing},  \cT_\prec} (\lambda) :=    dV_n \; d \Theta_{n-1} \; d\Omega_{n-1}  
\prod_{e = \{i,j\} \in E(  \cT_\prec )}    \big(  ( {\textbf  v}_i ^\e(\tau_e^-) - {\textbf  v}_j ^\e (\tau_e^-) ) \cdot\omega_{e }\big)_+ \,,
\end{align}
where the collision times are ordered according to the edges in the tree  (i.e.\,$d\Theta_{n-1}$ is supported on a simplex).
Iterating the change of variables \eqref{eq: identification measures} for the $|\lambda|-1$ clustering conditions prescribed by a given tree $ \cT_\prec$, one gets
\begin{equation}
\label{eq: changement variables  simplified Lambda}   
\begin{aligned}
\mu_\eps^{|\lambda|-1} \int d Z_n \; F^{(H)} (\lambda) \;  \indc_{ \lambda \ \hbox{\scriptsize min $\eps$-cluster path compatible with  } \cT_\prec } 
 = 
\int d y \; d\mu^\eps_{{\rm sing}, \cT_\prec} (\lambda) \; F^{(H)} (\lambda)
\; \indc_{\lambda \ \hbox{\scriptsize  min $\eps$-cluster path } }
 \, .
 \end{aligned}
\end{equation} 
Note that the parameters drawn from the measure $d\mu^\eps_{{\rm sing},  \cT_\prec}$ do not depend on $\eps$, thus for some values of $\eps$ the corresponding trajectory~${\textbf Z}^\eps_n([0,T]) $ may not form a min $\eps$-cluster path 
and the indicator function $\indc_{\lambda \hbox{ \scriptsize min $\eps$-cluster path } }$ in the right-hand side of \eqref{eq: changement variables  simplified Lambda} imposes this compatibility constraint.

\medskip As we now intend to take the limit $\mu_\eps \to \infty$ in the $\eps$-cluster paths, we emphasise the dependence on $\eps$ by writing  $\gl^\eps = \gl^\eps(n, \cT_\prec,y, V_n ,  \Theta_{n-1}, \Omega_{n-1})$ for a given ordered tree $ \cT_\prec$ and parameters~$(y, V_n ,  \Theta_{n-1}, \Omega_{n-1})$.
 Then $\gl^\eps$ converges, when $\mu_\eps$ tends to $\infty$, to a {\it limiting cluster path}~$\gl= \gl (n, \cT_\prec,y, V_n ,  \Theta_{n-1}, \Omega_{n-1})$  such that the  positions of the colliding  particles coincide at the collision times as in the definition below.
\begin{Def}[Limiting cluster path]
\label{Def: limiting trajectories}
Fix $n$, an ordered tree $ \cT_\prec$ of size $n$,  and  a  collection of parameters $(y , V_n,  \Theta_{n-1}, \Omega_{n-1})$
as in \eqref{eq: sing measure}.
The corresponding  \emph{limiting cluster path} $\lambda$ can be constructed as follows. 
In between two collision times, all particles evolve according to the free flow and 
for each edge $e = \{i,j\}  \in \cT_\prec$,  the corresponding constraints
are imposed at the collision time $\tau_e$:
\begin{itemize}
\item 
a collision occurs between the particles $i,j$, 
\item the positions of both particles coincide $x_i ( \tau_e ) = x_j  ( \tau_e )$,
\item   the velocities~$v_i ( \tau_e^+ )$ and $v_j  ( \tau_e^+ )$ are scattered  according to the  rule \eqref{eq: scattlaw} with scattering vector $ \omega_{e }$.
\end{itemize}

The associate singular measure is defined by
\begin{align}
\label{eq: sing measure limit}
 d\mu _{{\rm sing},  \cT_\prec} (\lambda) :=     dV_n \; d \Theta_{n-1} \; d\Omega_{n-1}  
\prod_{e = \{i,j\} \in E(  \cT_\prec )}   \big(  ( {\textbf  v}_i  (\tau_e^-) - {\textbf  v}_j  (\tau_e^-) ) \cdot\omega_{e }\big)_+ \,.
\end{align}
\end{Def} 
Recall  that the distance $\| \cdot \|_T$ is defined in  \eqref{eq: norm T}.
By construction, the following uniform convergence of the $\eps$-cluster paths $\lambda^\eps$ towards $\lambda$ holds
\begin{equation}
\label{eq: uniform convergence particles}
\| \lambda^\eps(n,\cT_\prec, y, V_n ,  \Theta_{n-1}, \Omega_{n-1}) - \lambda (n,\cT_\prec, y, V_n ,  \Theta_{n-1}, \Omega_{n-1})\|_T
\xrightarrow[\mu_\eps \to \infty]{} 0
\end{equation}
for all~$n$ and
$\cT_\prec$, and locally uniformly in~ $(y, V_n ,  \Theta_{n-1}, \Omega_{n-1})$  outside a set of measure 0.
Thanks to the continuity assumption~\eqref{eq: hypotheses h continuity} on the test function $H$, the limit of $d\mesh0 ( \gl ) $ in the Boltzmann-Grad limit  is the singular measure 
\begin{equation}
\label{forest-sing-measure}
 d{\pmb \nu}_T^{(H)} (\lambda) := \frac{1}{| \gl | \, !} \sum_{\cT_\prec \in \pmb{\cT}_{|\lambda|} ^{\prec}}  F^{(H)} (\gl)
\;  \;  d y \,d\mu _{{\rm sing},  \cT_\prec} (\lambda)
=:    dy \, d\tilde {\pmb \nu}_T^{(H)} (\lambda)\;.
\end{equation}
Note that in the limit there is no restriction on the set of parameters $(y, V_n,  \Theta_{n-1}, \Omega_{n-1})$ since by Proposition \ref{Prop: no cycle terms}, we know that   cycles have a vanishing probability.

\bigskip

We proceed exactly the same way with minimal $\eps$-aggregates $\{\lambda^\eps_1, \dots \lambda^\eps_k\}$. 
We recall that the overlaps occurring at time $0$ have been discarded in 
Proposition \ref{Prop: no cycle terms} so that all the clusterings occur only dynamically. 
 Assume that the overlaps between the $\eps$-cluster paths are prescribed by the
ordered tree $ \cT_\prec^{ov} $. For each edge $e=\{i, j\} $, let $q_i, q_j$ be the overlapping particles  in the cluster paths $ \lambda^\eps_i$ and $\lambda^\eps_j$. For fixed~$\eps>0$, the overlap  condition associated with the edge $e $   takes the form 
\begin{equation}
 \omega_{e } := \frac{ {\textbf  x}_{q_i }^\e(\tau_e) - {\textbf  x}_{q_j}^\e (\tau_e  )} {\eps}  \in {\mathbb S}^{d-1} \,,
\end{equation}
denoting by  $\tau_e$ the infimum of the  overlap  times between $ \lambda^\eps_i$ and $ \lambda^\eps_j$.
Note that the collision time used in \eqref{eq: omega collision} is uniquely defined, instead two particles overlap 
during a (short) time interval so that the overlapping time $\tau_e$ has to be more carefully prescribed  as an infimum.
Recall that before the overlap, the $\eps$-cluster paths~$\lambda^\eps_i,  \lambda^\eps_j$ are connected to two distinct components of the dynamical interaction graph which can move rigidly   with  the roots of the cluster paths  $y_i$ and~$y_j$.
Then we define the local change of variables
\begin{equation}
\label{eq: change variables2}
\hat y_e = y_j - y_i \in\bbT^d \mapsto(\tau_{e } ,  \omega_{e }) \in [0,T] \times  {\mathbb S}^{d-1}
\end{equation}
with inverse Jacobian determinant $ \mu_\eps^{-1}\big( ( {\textbf  v}_{q_i} ^\e(\tau_e^-) - {\textbf  v}_{q_j }^\e(\tau_e^-) )\cdot \omega_{e}\big)_+$. 
 Applying iteratively this change of variables for each edge of the  collision tree~$\cT_\prec^{ov}$, the  trajectories ${\textbf Z}^\eps_n ( [0,T])$ with $n =|\lambda^\eps_1|+\dots +|\lambda^\eps _k|$ can be built by the following parameters~:
\begin{itemize}
\item $y$  center of mass of the~$n$ positions at time~0 (also called root of the aggregate),
\item $(v_1, \dots, v_n)$ the particle velocities at the initial time,
\item for all $i \leq k$, an ordered tree $\cT_\prec^{(i)}$, and $(\Omega^{(i)}_{|\lambda_i|-1} ,\Theta_{|\lambda_i| -1}^{(i)}) $ parametrizing the collisions in the cluster paths $\lambda^\eps_i$,
\item an ordered tree $\cT_\prec^{ov}$, and $(\Omega_{k-1}, \Theta_{k-1})$ parametrizing the overlaps between the cluster paths.
\end{itemize}
Given all these parameters and outside a set of zero measure,  ${\textbf Z}^\eps_n  ( [0,T])$ converges uniformly in $[0,T]$ (in the sense of \eqref{eq: uniform convergence particles}) to a limiting configuration 
${\textbf Z}_n ( [0,T])$  such that 
\begin{itemize}
\item the  dynamics inside each limiting cluster path $(\lambda_i) _{i \leq k}$ is prescribed by Definition \ref{Def: limiting trajectories}, 
\item at an overlap time $\tau_e$ with $e = \{i,j\}$,  the overlapping particles  $q_i\in  \lambda_i, \ q_j\in \lambda_j$ coincide, i.e.  ${\textbf  x}_{q_i }( \tau_e ) = {\textbf  x}_{q_j}  ( \tau_e )$ (but their velocities are not scattered).
\end{itemize}

Then the limit of $d\mesh0 ( \gl_1 ) \dots d\mesh0 ( \gl_k ) \; \varphi_\eps^{\hbox{\rm \scriptsize min}} (\gl_1, \dots, \gl_k )  $ is the singular measure 
\begin{equation}
\label{cluster-sing-measure}
(-1)^{k-1} \sum_{\cT^{ov} _\prec}  dy\, d \Theta_{k-1} \; d\Omega_{k-1}  
\left( \prod_{ \{i,j\} \in E(  \cT_\prec ^{ov} )}    \big(  ( {\textbf  v}_{q_i} (\tau_e^-) - {\textbf  v}_{q_j } (\tau_e^-) ) \cdot\omega_{e }\big)_+ \right) d\tilde {\pmb\nu}_T^{(H)} (\lambda_1) \dots d\tilde {\pmb\nu}_T^{(H)} (\lambda_k)\,.
\end{equation}

\medskip

A limiting expression of the partition function is obtained below by 
combining \eqref{eq: simplified expansion} and the previous results.
\begin{Prop}
\label{prop: cluster expansion limit}
Let $T>0$ be the convergence time obtained in Proposition \ref{prop: cluster expansion}.
Then the following limit holds
\begin{equation}
\label{eq: cluster expansion - limit}
\lim_{\eps \to 0} \frac{1}{\mu_\eps} \log \cZ^\eps_T (e^H) 
=   \sum_{k \geq 1}\frac{(-1)^{k-1} }{k!}   \sum_{\cT^{ov} _\prec}  
\int dy \,d \Theta_{k-1} \; d\Omega_{k-1}  
\left( \prod_{ \{i,j\} \in E(  \cT_\prec ^{ov} )}    \big(  ( {\textbf  v}_i (\tau_e^-) - {\textbf  v}_j  (\tau_e^-) ) \cdot\omega_{e }\big)_+ \right) d\tilde {\pmb\nu}_T^{(H)} (\lambda_1) \dots d\tilde {\pmb\nu}_T^{(H)} (\lambda_k)  .
\end{equation}
Furthermore, as 
the series is uniformly absolutely convergent,  the limiting functional  $L (e^H)$  is such that $u \in \R \mapsto L(e^{uH}) $ is analytic in the neighborhood of 0, 
for $H$ satisfying \eqref{eq: hypotheses h}-\eqref{eq: hypotheses h continuity}.
\end{Prop}

Recalling \eqref{eq: log modified partition function bis}, the exponential moment 
$\gL^\eps_T (e^H) = \frac{1}{\mu_\eps} \left( \log \cZ^\eps_T (e^H) - \log \cZ^\eps_T (1) \right)$
can be controlled by the partition functions. 
From Proposition \ref{prop: cluster expansion limit}, its limiting expression is obtained 
as a series expansion in $e^{H}$.
In particular, if $H$ has the form (\ref{eq: hypotheses h bis}), this applies as well
to the exponential moment $\widetilde \gL^\eps_T (e^h)$  
defined in \eqref{eq: log partition function}.

\begin{Rmk} 
\label{Rmk: symetrisation} 
Following Remark \ref{global parametrization - rmk}, we can rewrite  $\int d\mesh0 ( \gl_1 ) \dots d\mesh0 ( \gl_k ) \; \varphi_\eps^{\hbox{\small min}} (\gl_1, \dots, \gl_k ) $ in \eqref{eq: simplified terms cluster expansion} as well as its limit with a more global parametrization by treating collisions and overlaps in a more symmetric way.

Ordering all the clustering conditions (coming both from collisions and overlaps) according to the forward flow, we index them by a single signed ordered tree $\overline {\cT_\prec}$ such that the $n$ particles form the vertices and each edge $e$ has a sign $s_e = +1$ if it records a collision or $s_e = -1$ if it records an overlap. 
Taking limits, we obtain 
\begin{equation}
\label{eq: simplified term Lambda}   
\lim_{\eps \to 0} \frac{1}{\mu_\eps}  \log \widetilde \cZ^\eps_T (e^h) =  
 \sum_{n =1}^\infty \frac{1 }{ n !} \sum_{\overline \cT_\prec \in \pmb{\cT}_n^{\prec, \pm}}
\int d y  \; d\mu_{sing, \overline \cT_\prec} (  {\textbf Z}_n )   \prod_{i \leq |\lambda|} f^0({\textbf z}_i(0)) \exp \Big( h   \big( {\pmb z}_i ([0,T]) \big)  \Big),
\end{equation}
where $\pmb{\cT}_n^{\prec, \pm}$ is the set of signed ordered   trees with $n$ vertices, 
 and 
\begin{equation} \label{eq: sing measure limit with signs}
d\mu _{{\rm sing},  \overline \cT_\prec} (  {\textbf Z}_n ) :=     dV_n \; d \Theta_{n-1} \; d\Omega_{n-1}  
\prod_{e = \{i,j\} \in E(  \overline \cT_\prec )}  s_e  \big(  ( {\textbf  v}_i  (\tau_e^-) - {\textbf  v}_j  (\tau_e^-) ) \cdot\omega_{e }\big)_+ \,.
\end{equation}
\end{Rmk}

Recall that 
$\widetilde \gL^\eps_T (e^h) 
 = \frac{1}{\mu_\eps} \left( \log \widetilde\cZ^\eps_T (e^h) - \log \widetilde\cZ^\eps_T (1) \right)$.
Then, as a direct consequence of \eqref{eq: simplified expansion}, we  recover Theorem 6 of \cite{BGSS2} which is stated below.
\begin{Cor}
\label{prop: cluster expansion limit h}
Let $T>0$ be the convergence time obtained in Proposition \ref{prop: cluster expansion}.
Then the following limit holds
\begin{equation}
\label{eq: limite moment exponential h}
\lim_{\eps \to 0} \widetilde \gL^\eps_T (e^h)  
= -1 +  \sum_{n =1}^\infty \frac{ 1 }{ n !} 
\sum_{\overline \cT_\prec \in \pmb{\cT}_n^{\prec, \pm}}
\int dy\, d\mu_{{\rm sing}, \overline \cT_\prec} \left(  {\textbf Z}_n  \right) 
   \;  \prod_{i = 1}^n f^0\big( {\textbf z}_i (0) \big) \exp \Big(  h  \big( {\textbf z}_i ( [0,T]) \big) \Big)  .
\end{equation}
Furthermore the limiting functional $\widetilde \Lambda_T (e^h)$    is such that $u \in \R \mapsto \widetilde \Lambda_T (e^{uh}) $ is analytic in the neighborhood of 0, 
for $h$ satisfying~\eqref{eq: hypotheses h bis}.
\end{Cor}


\begin{Rmk}
The series expansion \eqref{eq: limite moment exponential h} is a key tool in 
{\rm\cite{BGSS2}} to prove that the functional $t \in [0,T] \mapsto \widetilde \gL^\eps_t$ can be characterised (on a nice class of test functions) in terms of a Hamilton-Jacobi equation.
This Hamilton-Jacobi equation provides then a refined information on the particle system : in particular, it encodes
the fluctuating Boltzmann equation and the large deviations (quantifying atypical particle evolutions).
We refer to  {\rm\cite{BGSS2}} for the mathematical details 
and to {\rm\cite{BGSS1, BGSS_ICM}} for an overview of these results and a discussion of their physical interpretation.
\end{Rmk}



\section{Typical behaviour of the density and of the particle trajectories}
\label{sec: limiting equations}

\subsection{Derivation of the Boltzmann equation}
\label{sec: Boltzmann equation}

As a first application of Proposition \ref{prop: cluster expansion limit}, we are going to recover 
Theorem \ref{thm: Lanford}, i.e.\,that the limiting density of the hard sphere dynamics follows the Boltzmann equation for short times.  This result is restated below in Proposition \ref{Prop: convergence Boltzmann}, making the  link with our previous discussions.
We start by introducing a notion of strong solution of  the Boltzmann equation \eqref{eq:Beq}. For simplicity, 
we will use a shorthand notation for the collision operator  and rewrite  \eqref{eq:Beq} as 
\begin{equation}
\partial_t f +v \cdot \nabla_x f = \pmb{C}(f,f).
\label{eq:Beq collision C}
\end{equation}
By a fixed point argument   \cite{Ukai1974,Ukai1976, Kaniel_Shinbrot}, one can show  that under the assumptions \eqref{eq: Gaussian bound} on $f^0$ there exists a unique, stable solution of  the Boltzmann equation on a time interval $[0,T^\star]$ with~$ T^\star \geq T$, the convergence time obtained in Proposition \ref{prop: cluster expansion}.  
In particular, this solution is a mild solution and  takes the following form for $t \leq T^\star$
\begin{align}
\label{eq: mild solution}
f (t)  & = \pmb{S}_1(t) f^0  + \int_0^t d \tau \, \pmb{S}_1 ( t-  \tau) \, \pmb{C} \big( f (\tau), f( \tau ) \big) ,
\end{align}
where the operator $\pmb{S}_1 (u)$ acts as the backward free transport during time $u$.

\medskip

\begin{Prop}
\label{Prop: convergence Boltzmann}
Let $T>0$ be the convergence time obtained in Proposition \ref{prop: cluster expansion}. The empirical measure converges to the solution
of the Boltzmann equation in the following sense : for any $t \leq T$,  any continuous  test function 
$h$ in $L^\infty(\bbT^d \times \bbR^d)$ and~$\gd >0$
\begin{equation}
\label{eq: convergence Boltzmann en proba}
\bbP_\eps \left( \Big|\pi^\eps_t(h) -  
\int_{\bbT ^d\times \bbR^d} dz \, f(t,z) h(z)\Big| > \delta \right)
\xrightarrow[\mu_\eps \to \infty]{} 0 \,.
\end{equation}
\end{Prop}
	
\medskip

\noindent
{\it Proof.}
To prove Proposition \ref{Prop: convergence Boltzmann}, it is enough to show the convergence in law, i.e. the limit
\begin{equation}
\label{eq: density at time t}
\lim_{\eps \to 0} \bbE_\eps \left[  \pi^\eps_t (h)   \right] 
= \int_{\bbT ^d\times \bbR^d} dz \, f(t,z) \, h(z).
\end{equation}
Indeed the convergence in probability \eqref{eq: convergence Boltzmann en proba} follows from the Markov inequality and the $L^2$ inequality  
\begin{equation}
\label{eq: covariance moyenne Boltzmann}
\bbE_\eps \left[ \left(  \pi^\eps_t (h) - \int dz \, f(t,z) \, h(z)  \right)^2 \right]
\leq  O \left( \frac{\| h \|_\infty^2 }{\mu_\eps}  \right) +
2 \left( \bbE_\eps \left[  \pi^\eps_t (h)   \right] - \int_{\bbT ^d\times \bbR^d} dz \, f(t,z) \, h(z) \right)^2 ,
\end{equation}
which  is a consequence of  the concentration estimate \eqref{eq: covariance pi}.

We turn now to the derivation of \eqref{eq: density at time t}.
By Propositions \ref{prop: cluster expansion} and \ref{prop: cluster expansion limit}, the functionals
$\widetilde \gL^\eps_T$ and its limit $\widetilde  \gL_T$ are analytic so that 
\begin{equation}
\label{eq: derivation limite}
\lim_{\eps \to 0} \bbE_\eps \left[  \pi^\eps_{[0,T]}(h)   \right] 
= \lim_{\eps \to 0} \partial_u \, \gL^\eps_T (e^{u \, h}) \Big|_{u =0}
=  \partial_u \, \gL_T (e^{u \, h}) \Big|_{u =0} .
\end{equation}
For a given $t \in [0,T]$, we choose   
$$
h  \big( {\textbf z} ( [0,T] ) \big) = h  \big( {\textbf z} ( t) \big),
$$
with a continuous  test function $h$ in $L^\infty( \bbT^d \times \bbR^d)$.
Note that the identity \eqref{eq: derivation limite} is obtained by derivation, so that instead of  conditions \eqref{eq: hypotheses h}-\eqref{eq: hypotheses h continuity}, it is enough to assume that the test function $h$ is continuous and bounded.
In the previous sections, the time window $[0,T]$ was fixed, but it will be convenient 
to reduce it now to $[0,t]$.
We deduce, from the identity \eqref{eq: derivation limite}, the limiting counterpart of \eqref{eq: average Pi}
(using the notation of \eqref{eq: simplified term Lambda})
\begin{align}
\lim_{\eps \to 0} \bbE_\eps \left[  \pi^\eps_t (h)   \right] 
& =
\sum_{n =1}^\infty \frac{ 1 }{ n !} 
\sum_{\overline \cT_\prec \in \pmb{\cT}_n^{\prec, \pm}}
\int dy \; d \mu_{{\rm sing}, \overline \cT_\prec}^{[0,t]} \left(  {\textbf Z}_n \right)
\; \left( \sum_{i =1}^n h  \big( {\textbf z}_i (t) \big) \right)   
\;  \prod_{i = 1}^n f^0\big( {\textbf z}_i (0) \big) \nonumber \\
& =
\sum_{n =1}^\infty \frac{ 1 }{ (n-1) !} 
\sum_{\overline \cT_\prec \in \pmb{\cT}_n^{\prec, \pm}}
\int   dx_1 \; d \mu_{{\rm sing}, \overline \cT_\prec}^{[0,t]} \left( {\textbf Z}_n  \right)
\;  h  \big( {\textbf z}_1 (t) \big)   \;  \prod_{i = 1}^n f^0\big( {\textbf z}_i (0) \big),    
\label{eq: limite moment 1}
\end{align}
where the singular measure $\mu_{{\rm sing}, \overline \cT_\prec}^{[0,t]}$ is defined as in \eqref{eq: sing measure limit with signs}, with the superscript indicating that the limiting cluster paths ${\textbf Z}_n={\textbf Z}_n\left([0,t]\right)$  are here restricted to the time interval $[0,t]$.
The second equality is obtained by the symmetry of the particles and 
 by using the initial position $x_1$ as a root instead of the center of mass $y$.

\medskip

To complete \eqref{eq: density at time t}, it remains to show that 
the limiting particle distribution at time $t$ defined by the right-hand side of \eqref{eq: limite moment 1}
coincides with the solution $f(t)$ of the Boltzmann equation.
We proceed in 2 steps : 
first the series  \eqref{eq: limite moment 1} is simplified as many terms cancel out, 
leading  to a representation of the particle distribution 
which is then shown to satisfy the Boltzmann  equation~\eqref{eq: mild solution}.

\medskip

\noindent
\textbf{Step 1} : {\it Restriction to a relevant time ordered cluster of influence.}
 
\begin{figure}[h] 
\centering
\includegraphics[width=3.2in]{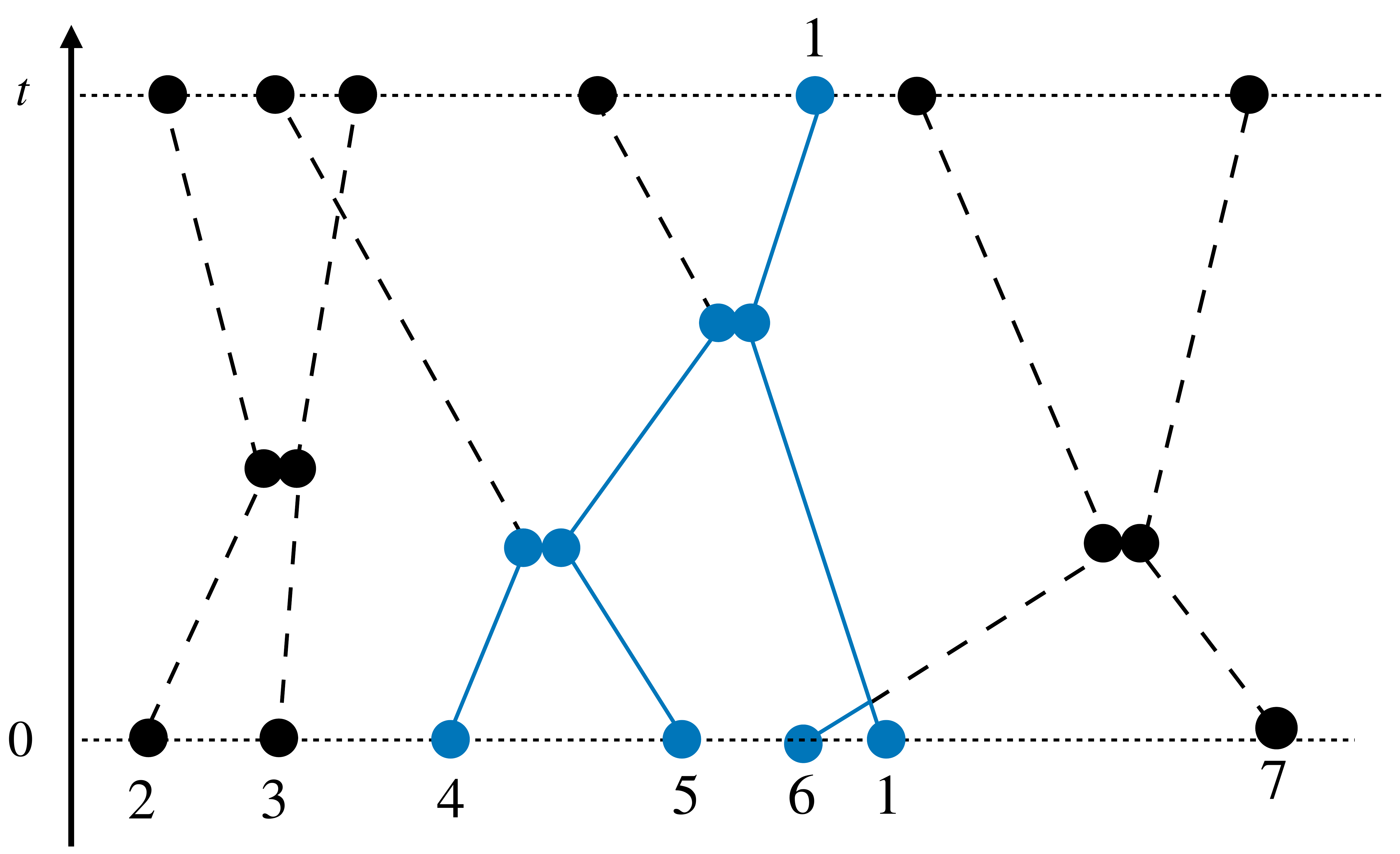} 
\hskip1cm
\includegraphics[width=3.2in]{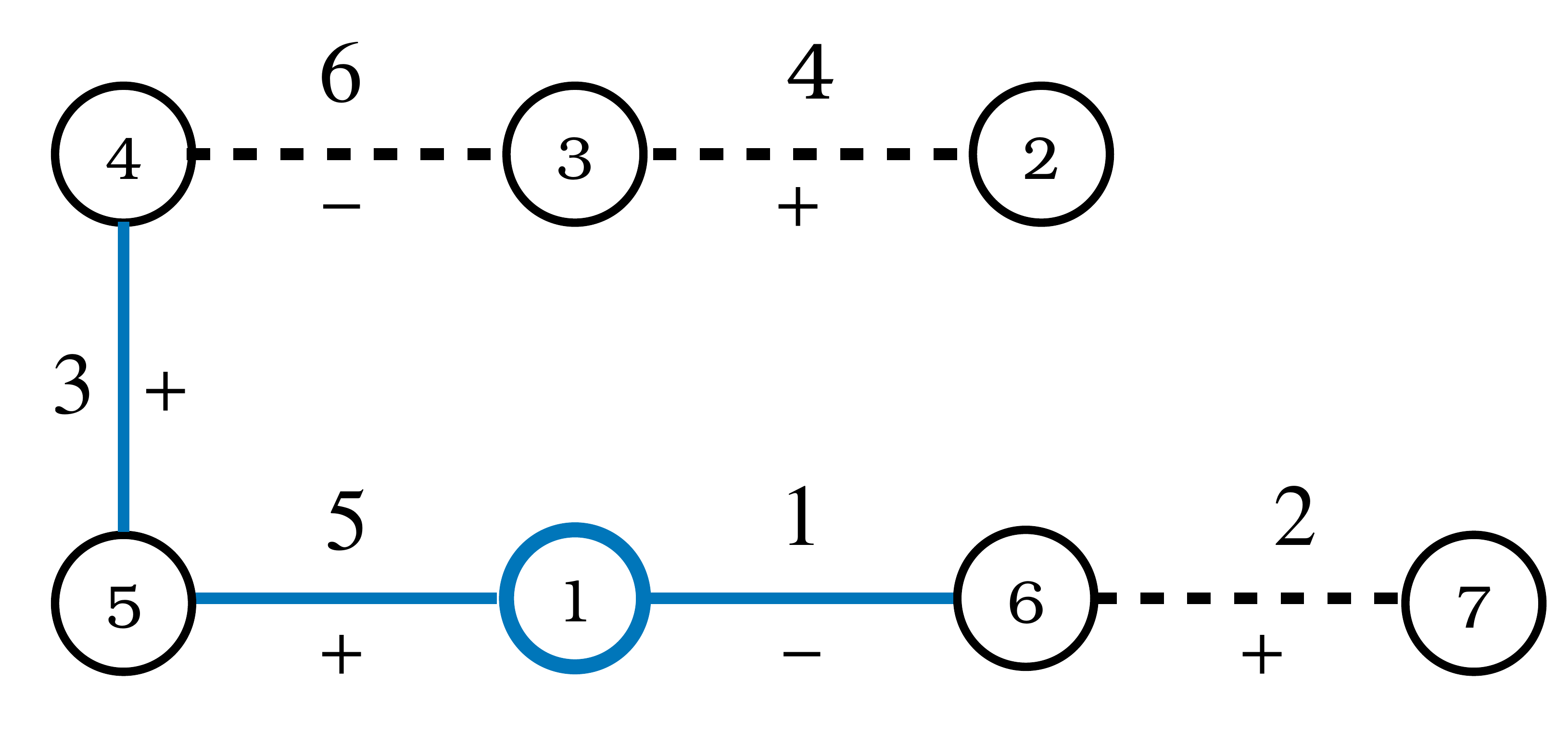} 
\caption{\small On the left, particle trajectories associated with a term of order $n =7$ in \eqref{eq: limite moment 1} are depicted and
the corresponding ordered graph $\overline \cT_\prec$ is depicted on the right : each edge has an order and a sign $\pm$ to record if it is a  collision or an overlap. The relevant part $\cA$ of the trajectories, which determine ${\textbf z}_1(t)$, is represented by blue lines and the edges which can be neglected by dashed lines.
Ultimately the initial coordinates of the particles in the blue tree $\cA$ and the sign of the clusterings prescribe ${\textbf z}_1 (t)$.
}
\label{figure: arbre reduction}
\end{figure}

First, we are going to simplify each term  of the series  \eqref{eq: limite moment 1} by showing that only the particles in a cluster of influence of particle 1 are needed to compute $h  \big({\textbf z}_1 (t) \big)$ (see Figure \ref{figure: arbre reduction}).
To extract the relevant information, we consider a signed ordered tree $\overline \cT_\prec$ of size $n$ and proceed  by building recursively a growing collection of subtrees $\cA_1 \subset \cA_2 \subset \dots$ as follows.
Starting from the vertex $\cA_1 = \{ 1 \}$ associated with the particle $1$, all the vertices and edges connected to $1$ are added to form the set $\cA_2$. Suppose that $\ell$ belongs to $\cA_2$ and that the edge 
$\{1,\ell\}$ has order $k$ 
then all the neighbours of $\ell$ are added to $\cA_3$
provided they are linked to $\ell$ by an edge with order smaller than $k$, i.e.\,if the corresponding clustering has occurred before $k$. Iterating this procedure leads to an ordered tree $\cA$ from which the configuration ${\textbf z}_1 (t)$ can be recovered. We will denote by $\pmb{\cA}_n^{\prec, \pm}$ the set containing the trees with $n$ vertices of the previous form, i.e. the trees rooted in $1$ such that the edge orders are decreasing when examined from the root to a leaf (see Figure \ref{figure: arbre reduction}).

Given $\cA \in \pmb{\cA}_n^{\prec, \pm}$, let $\overline \cT_\prec$ be a tree in which $\cA$ can be embedded but which has at least one more edge than $\cA$. Choosing one leaf in  $\overline \cT_\prec$ which is not in $\cA$, one can build another tree $\overline \cT_\prec'$ by simply changing the sign of the edge connecting this leaf. This changes the tree locally without influencing the value of $h  \big( {\textbf z}_1 (t) \big)$ so that
\begin{align}
\label{eq: compensation}
\int   dx_1 \; d \mu_{{\rm sing}, \overline \cT_\prec'}^{[0,t]} \left(  {\textbf Z}_n\right)
\;  h  \big( {\textbf z}_1 (t) \big)   \;  \prod_{i = 1}^n f^0\big( {\textbf z}_i (0) \big)
= - \int   dx_1 \; d \mu_{{\rm sing}, \overline \cT_\prec}^{[0,t]} \left(  {\textbf Z}_n \right)
\;  h  \big( {\textbf z}_1 (t) \big)   \;  \prod_{i = 1}^n f^0\big( {\textbf z}_i (0) \big).  
\end{align}
Thus the sum \eqref{eq: limite moment 1} reduces to 
\begin{align*}
\lim_{\eps \to 0} \bbE_\eps \left[  \pi^\eps_t (h)   \right] 
=
\sum_{n =1}^\infty \frac{ 1 }{ (n-1) !} 
\sum_{\cA \in \pmb{\cA}_n^{\prec, \pm}}
\int   dx_1 \; d \mu_{{\rm sing}, \cA}^{[0,t]} \left(  {\textbf Z}_n \right)
\;  h  \big( {\textbf z}_1 (t) \big)   \;  \prod_{i = 1}^n f^0\big( {\textbf z}_i (0) \big).   
\end{align*}
Note that the time ordering of the edges is sufficient to recover from $\cA$ the final position ${\textbf z}_1 (t)$. 
Relabelling the vertices according to the backward order of clusterings,  and denoting by~$\pmb{ \hat A}_n^{\prec, \pm}$ the set of rooted trees with ordered and signed edges,  we finally obtain 
\begin{align}
\lim_{\eps \to 0} \bbE_\eps \left[  \pi^\eps_t (h)   \right] 
=
\sum_{n =1}^\infty \sum_{\hat \cA \in \pmb{ \hat A}_n^{\prec, \pm}}
\int   dx_1 \; d \mu_{{\rm sing}, \hat \cA}^{[0,t]} \left( {\textbf Z}_n \right)
\;  h  \big( {\textbf z}_1 (t) \big)   \;  \prod_{i = 1}^n f^0\big( {\textbf z}_i (0) \big).
\label{eq: limite simple tree}
\end{align}
This defines, by duality, $g(t,z)$ the typical distribution  of a particle at time $t$
by imposing a constraint on the final coordinates 
\begin{align}
g(t,z) =
\sum_{n =1}^\infty \sum_{\hat \cA \in \pmb{ \hat A}_n^{\prec, \pm}}
\int d \mu_{{\rm sing}, \hat \cA}^{[0,t]} \left( {\textbf Z}_n \right)
\;  \delta_{z -  {\textbf z}_1 (t)}   \;  \prod_{i = 1}^n f^0\big( {\textbf z}_i (0) \big),
\label{eq: limite distribution Boltzmann def}
\end{align}
 where the degree of freedom of the initial position $x_1$ is fixed by the constraint at time $t$.

To give a meaning to the singularity in the collision operator of the Boltzmann equation \eqref{eq: mild solution},  
 some  regularity properties on $g$ will be needed.
For this we suppose, for a moment, that the initial distribution $f^0$ satisfies the following additional smoothness assumption
for any $\ell \leq d + 1$
\begin{equation}
\label{eq: smoothness Gaussian bound} 
\forall z \in \bbT ^d\times\bbR^d, \qquad 
|\partial_{x_{i_1}} \dots \partial_{x_{i_\ell}}  f^0(z)| \leq C \cM_{\beta/2} (v).
\end{equation}
We can then deduce that $g$ is continuous. 
\begin{Lem}
\label{Lem: continuity}
Assuming  \eqref{eq: smoothness Gaussian bound}, then  the density $(x,v) \mapsto g\big( t,(x,v) \big)$ belongs to $C^0_x \big( {\textbf M}^1_v \big)$ for all $t \in [0,T]$, where $ {\textbf M}^1_v$ is the space of measures such that $(1+|v|)$ is integrable.
\end{Lem}
The derivation of this lemma is postponed until the end of the proof of Proposition \ref{Prop: convergence Boltzmann}. Note that assumption \eqref{eq: smoothness Gaussian bound} will be lifted at the end of the proof by a density argument.

\medskip

\noindent
\textbf{Step 2} : {\it Identification of the Boltzmann equation.}

We are going to check that $g$ is a mild solution of the form \eqref{eq: mild solution}.
For this, it is convenient to interpret the series \eqref{eq: limite distribution Boltzmann def} defining~$g(t,z)$
as a backward evolution. This point of view coincides with the standard notion of pseudo-trajectories 
 used to derive Lanford's Theorem (see e.g.\,\cite{Cercignani_Illner_Pulvirenti}).
Fixing at time $t$ the position $z = {\textbf z}_1(t)$ of particle $1$, the collision tree $\hat \cA$ is built backward by adding  the particles which are interacting dynamically with particle $1$ (cf.\,Figure \ref{figure: arbre reduction}). 
Proceeding from time $t$, particle $1$  follows the backward free flow either up to time 0 or up to a time 
 $\tau$ at which a branching occurs, say with particle 2. The corresponding edge 
in the tree $\hat \cA$ is associated with a sign $-$ if this event is an overlap
and $+$ if it is a collision (in which case scattering occurs). 
Removing this edge splits $\hat \cA$ into two smaller trees $\hat \cA_1, \hat \cA_2$, containing respectively $1$ and $2$,
 which encode the rest of the trajectory on $[0,\tau]$.
 Thus the singular measure on $[0,t]$ is the product of the singular measures of the two trees 
 in $[0,\tau]$ with the constraint ${\textbf  x}_1 (\tau) = {\textbf  x}_2 (\tau)$ 
\begin{align*}
\mu_{{\rm sing}, \hat\cA}^{[0,t]} \left(  {\textbf Z}_n \right)
= \mu_{{\rm sing}, \hat\cA_1}^{[0,\tau]} \left(  {\textbf Z}_{n_1} \right)
\mu_{{\rm sing}, \hat\cA_2}^{[0,\tau]} \left(  {\textbf Z}_{n_2} \right)
\; s \big(   ( {\textbf  v}_1(\tau^+) - {\textbf  v}_2 (\tau^+) )\cdot \omega \big)_+
\delta_{{\textbf  x}_1(\tau) - {\textbf  x}_2(\tau)} \; d\tau \; d\omega  ,
\label{eq: splitting measure}
\end{align*}
where $\omega$ is the scattering parameter for the last encounter at time $\tau$ which can be a collision or an overlap according to $s = \pm1$. We therefore distinguish the velocities at times  $\tau^+$ and $\tau^-$.
In this way, \eqref{eq: limite distribution Boltzmann def} can be rewritten
\begin{align*}
g(t,z)  = f^0\big( x-v t, v\big)  + 
\sum_{s = \pm 1}
\sum_{n_1,n_2 =1}^\infty \sum_{\hat \cA_1 \in \pmb{ \hat A}_{n_1}^{\prec, \pm} \atop
\hat \cA_2 \in \pmb{ \hat A}_{n_2}^{\prec, \pm}}
\int d\tau \; d\omega &
\int d \mu_{{\rm sing}, \hat\cA_1}^{[0,\tau]} \left(  {\textbf Z}_{n_1} \right)
\; d \mu_{{\rm sing}, \hat\cA_2}^{[0,\tau]} \left(  {\textbf Z}_{n_2} \right)
\;  \delta_{(x-v (t-\tau), v)  -  {\textbf z}_1 (\tau)}  \; \delta_{{\textbf  x}_1(\tau) - {\textbf  x}_2(\tau)} \\
&  \times    \; s \big(   ( {\textbf  v}_1(\tau^+) - {\textbf  v}_2 (\tau^+) )\cdot \omega \big)_+ \;  \prod_{i = 1}^{n_1+n_2} f^0\big( {\textbf z}_i (0) \big).
\end{align*}
In the integral, the free transport operator during time $t - \tau$ can be identified. 
Depending on $s= \pm1$, the 
velocities $(v,v_2)$ at time~$\tau^+$ may change at $\tau^-$. Thus we define 
\begin{equation*}
(v^s, v_2^s)  
= \begin{cases}
(v', v_2'), \qquad & \text{if $s=1$ (scattering occurs at the collision with deflection angle $\omega$)},\\
(v, v_2),  \qquad  & \text{if $s=-1$ (the velocities are unchanged by the overlap)}.
\end{cases}
\end{equation*}

Using the uniform convergence of the series, Fubini's theorem applies
\begin{align*}
g(t,z)  = f^0\big( x-v t, v\big)  
+ \sum_{s = \pm 1}  \int & d\tau \; d \omega \; d v_2 \; s \big(   ( {  v}_1 - {  v}_2  )\cdot \omega \big)_+\\
&  \times
\left( 
\sum_{n_1=1}^\infty \sum_{\hat \cA_1 \in \pmb{ \hat A}_{n_1}^{\prec, \pm} }
 \int d \mu_{{\rm sing}, \hat\cA_1}^{[0,\tau]} \left(  {\textbf Z}_{n_1} \right)
\;  \delta_{(x-v (t-\tau), v^s)  -  {\textbf z}_1 (\tau^-)} 
\;  \prod_{i = 1}^{n_1} f^0\big( {\textbf z}_i (0) \big) \right)\\
&  \times
\left( \sum_{n_2 =1}^\infty \sum_{\hat \cA_2 \in \pmb{ \hat A}_{n_2}^{\prec, \pm}}
\int d \mu_{{\rm sing}, \hat\cA_2}^{[0,\tau]} \left(  {\textbf Z}_{n_2} \right)
\;  \delta_{(x-v (t-\tau), v_2^s)  -  {\textbf z}_2 (\tau^-)}  \; 
\;  \prod_{i = 1}^{n_2} f^0\big( {\textbf z}_i (0) \big) \right)
\end{align*}
 and the product of the densities $g(\tau^-)$ can be recovered. 
The densities are independent as they are defined by different sets of parameters.
The continuity, derived in Lemma \ref{Lem: continuity}, allows us to make sense of the Dirac 
condition on the positions at time~$\tau$.
The gain and loss terms in the collision operator arise from the parameter $s= \pm$. 
Thus $g$ solves the mild form \eqref{eq: mild solution} of the Boltzmann equation.

\medskip

To conclude the proof of Proposition \ref{Prop: convergence Boltzmann}, it remains to show 
that the regularity assumption \eqref{eq: smoothness Gaussian bound} on the initial data is not necessary.
Indeed the limiting density $g$ introduced in \eqref{eq: limite distribution Boltzmann def}  remains well defined as a measure which can be approximated by a sequence $(g^\delta)_{\delta>0}$ obtained by regularising $f^0$ by smooth densities $f^{0,\delta}$ satisfying~(\ref{eq: smoothness Gaussian bound}). By the previous argument each $g^\delta$ is a mild solution of the Boltzmann equation and the stability of the Boltzmann equation 
with respect to uniform convergence implies that the limit $g$ solves also \eqref{eq: mild solution}.
This completes the derivation of Proposition \ref{Prop: convergence Boltzmann}.

\bigskip

\noindent
{\it Proof of Lemma \ref{Lem: continuity}.}

First note that the integrability condition can be deduced from the assumption \eqref{eq: hypotheses h} as the test functions are allowed to diverge as $v^2$. 
To prove the statement on the continuity with respect to~$x$, let us first recall that the measure $\mu_{{\rm sing}, \hat \cA}^{[0,t]}$ prescribes a set of trajectories~${\textbf Z}_n ([0,t])$ which can be moved rigidly by shifting the root $y$ without changing the weights in the measure. Thus the final condition ${\textbf x}_1 (t) = x$ can be rewritten as a condition on the  initial position $x_1 = x + Y({\textbf Z}_n ([0,t]))$ where the last term (which will not be computed here) depends only on the internal structure of the rigid cluster.
For any $r \in \bbT^d$, the spatial shift operator~$\tau_r$ applied to the test function $h$ is defined by 
$\tau_r h(x) = h(x-r)$.
A shift by $r$ of the test function can be interpreted by duality by shifting rigidly by $-r$ the whole trajectory ${\textbf Z}_n ([0,t])$  so that
\begin{equation}
\label{eq: first shift}
\int dx dv \,  g \big( t,(x,v) \big) \,  \tau_r h (x,v) 
= 
\sum_{n =1}^\infty \sum_{\hat \cA \in \pmb{ \hat A}_n^{\prec, \pm}}
\int  d x_1 \; d \mu_{{\rm sing}, \hat \cA}^{[0,t]} \left( {\textbf Z}_n \right)
\;  h  \big( {\textbf z}_1 (t) \big)   \;   \prod_{i = 1}^n f^0\big( {\textbf z}_i (0) +(r,0) \big)
\end{equation}
where the  shift  is only felt at the level of the initial measure.
This boils down to applying a shift $\tau_{-r}$ to the function
\begin{equation}
\label{eq: shifted initial data}
x \mapsto  \Gamma_{{\textbf Z}_n(0)} (x) :=  \prod_{i = 1}^n f^0\big( {\textbf z}_i (0) + (x,0) \big)
\end{equation}
whose  derivatives of order $d+1$ are bounded thanks to \eqref{eq: smoothness Gaussian bound}
with a growth as $n^{d+1}$.
Thus, using the discrete derivative $\tau_r - \text{Id}$,  we deduce that uniformly in $r$ 
\begin{equation}
\left| \int dx dv \,  g \big( t,(x,v) \big) \;  \frac{1}{r^{d+1}} ( \tau_r - \text{Id})^{d+1}  h (x,v)  \right|
\leq C \; \Big\| \frac{h}{1+ |v|} \Big\|_\infty .
\end{equation}
This implies that~$d+1$ space derivatives  of~$g$ are measures in~$x$ with values in~$ {\textbf M}^1_v$, hence~$g$ is continuous   in~$x$ with values in $ {\textbf M}^1_v$.

\subsection{The dynamical cluster process}
\label{sec: tanaka}

The cluster paths play a key role in the cluster expansion of the exponential moment and it is interesting 
to study the dynamics of these cluster paths in its own right as a relevant observable of the hard sphere dynamics.
The methods developed in this paper provide a direct way for studying the particle trajectories as 
the exponential moments encode the statistics of the cluster paths.
In particular, an analogous result  to Proposition \ref{Prop: convergence Boltzmann} can be derived at the level of clusters. Before stating it, let us fix some notations.

\medskip

Following Definition \ref{Def: limiting trajectories}, a limiting cluster path of size $n$ is determined by a decorated cluster (tree) graph
\begin{equation}
\label{eq: cluster path}
\left(n,  \cT_\prec,  y , V_n,  \Theta_{n-1}, \Omega_{n-1}\right) 
\in  \TT  := \bigcup_{n \in \bbN} \{ n\} \times \pmb{\cT}_n ^{\prec} \times \bbT^d \times
\bbR^{d n} \times ( 0,\infty )_s^{n-1} \times \bbS^{(d-1) (n - 1)} ,
\end{equation}
where $( 0, \infty)_s^{n - 1}$ is a shorthand notation for the simplex 
$0< \tau_1< \tau_2 < \dots < \tau_{n -1} $.
 By abuse of notation, we shall use in this section the symbol $\gl$ for elements in $\TT$.
Notice that there is no restriction on time in \eqref{eq: cluster path}, nevertheless it will be convenient to define 
the notion of cluster paths \emph{formed} in the time interval $[0,t]$ by the additional constraint $\tau_{n -1} \leq t$.
We further recall that, for limiting cluster paths on $[0,T]$, \eqref{eq: cluster path} provides a convenient way to parametrise 
the test functions $H$ defined at the level of \eqref{eq: hypotheses h}; we shall write below $H = H(\gl), \gl \in \TT$ for such parametrization.

Consider two  limiting cluster paths formed in the time interval $[0,t]$, with decorated graphs
\begin{equation*}
\gl_1 = (n_1,  \cT_{\prec,1},  y_1 , V_{n_1},  \Theta_{ n_1 -1}, \Omega_{ n_1 -1}) 
\quad \text{and} \quad 
\gl_2 = (n_2,  \cT_{\prec,2},  y_2 , V_{n_2},  \Theta_{ n_2 -1}, \Omega_{ n_2 -1}) 
\end{equation*}
(such that $\Theta_{ n_1 -1}, \Theta_{ n_2 -1}$ are  restricted to $(0,t)$). 
By definition, a collision at time $t$ with deflection angle $\omega$ between particles $i \leq |\gl_1|$ and $j \leq |\gl_2|$ 
gives rise to the \emph{merging} of $\gl_1,\gl_2$   creating the limiting cluster path  with decorated graph
\begin{equation}
\label{eq: merging cluster path}
[\gl_1 \wedge \gl_2]^{i,j,t, \omega} 
:= \big(n_1 +n_2,  \cT_{\prec}^{ i,j},  y , V_{n_1+n_2},  \Theta_{ n_1 + n_2 -1}, \Omega_{ n_1+ n_2 -1} \big)
\end{equation}
 built as follows :
\begin{itemize}
\item $\Theta_{ n_1 + n_2 -1}$ contains the reordered collision times of $\gl_1, \gl_2$ as well as the time $t$ of the new collision,
\item the tree $\cT_{\prec}^{i,j}$ is the aggregation of the trees $\cT_{\prec,1}$ and $\cT_{\prec,2}$ obtained by creating a new edge $\{i,j\}$. The edges of $\cT_{\prec}^{i,j}$ are ordered according to the collision times in 
$ \Theta_{ n_1 + n_2 -1}$,
\item the root $y$ corresponds to the center of mass of the initial positions of the cluster 
$[\gl_1 \wedge \gl_2]^{i,j,t, \omega}$, and the initial (reordered) velocities are indicated by $V_{n_1+n_2}$,
\item the velocities of particles $i,j$ are scattered at time $t$  and the deflection angle $\omega$ associated with the new edge $\{i,j\}$ is added to the previous collection of deflection angles to create $\Omega_{ n_1+ n_2 -1}$.
\end{itemize}
By construction, the new cluster path is symmetric with respect to $\gl_1,\gl_2$.

\begin{thm}
\label{Prop: convergence dynamical clusters}
Let $T>0$ be the convergence time obtained in Proposition \ref{prop: cluster expansion}.
Then for any $t \leq T$, there is a measure $\tanaka{t}$  on $\TT$ (defined explicitly in 
\eqref{eq: def tanaka}) such that the empirical measure on the cluster paths converges to $\tanaka{t}$
 in the following sense :  for any test function $H$ satisfying \eqref{eq: hypotheses h}-\eqref{eq: hypotheses h continuity} and $\gd >0$
\begin{equation}
\label{eq: convergence clusters en proba}
\bbP_\eps \left( \Big| \Pi^\eps_t(H) -  \int_\TT  d \tanaka{t} (\gl_1) H(\gl_1) \Big| > \delta \right)
\xrightarrow[\mu_\eps \to \infty]{} 0 ,
\end{equation}
where the test functions  on the limiting cluster paths can be parametrised
as in \eqref{eq: cluster path}.

Furthermore,  under the regularity assumption \eqref{eq: smoothness Gaussian bound}  for $f^0$,  the limiting measure obeys an evolution equation which is the counterpart 
of \eqref{eq: mild solution} at the level of cluster paths
\begin{align}
\label{eq: mild tanaka}
 \forall t \leq T, \quad 
& \int  d \tanaka{t} (\gl_1) H (\gl_1) 
 =  \int  d \tanaka{0} (\gl_1) H (\gl_1)   \\
& \quad + \frac{1}{2} \int  d \tau  \, d \omega \, d \tanaka{\tau} (\gl_1) \, d \tanaka{\tau} (\gl_2) 
\sum_{ i \in \gl_1 \atop j \in \gl_2 }  \Big(   H \big( [\gl_1 \wedge \gl_2]^{i,j,\tau, \omega} \big) -   H (\gl_1) -   H (\gl_2) \Big) 
\, \delta_{{\textbf  x}_i (\tau) - {\textbf  x}_j(\tau)}   \,   \big(  ( {\textbf  v}_i (\tau) - {\textbf  v}_j  (\tau) ) \cdot \omega \big)_+ , \nonumber
\end{align}
where the clusters appearing in the argument of $H$ are formed in the time interval $[0,\tau]$ and merged at time $\tau \in [0,t]$.
\end{thm}

The singularity of the Dirac mass in the integral \eqref{eq: mild tanaka} makes sense thanks to the continuity of measures with respect to their roots which can be derived as in Lemma \ref{Lem: continuity}. 
We expect   the coagulation equation \eqref{eq: mild tanaka} to be stable and that 
Theorem \ref{Prop: convergence dynamical clusters}  holds  without the smoothness assumption \eqref{eq: smoothness Gaussian bound}.

Notice that the time-zero term in \eqref{eq: mild tanaka} is given by $$\int  d \tanaka{0} (\gl_1) H (\gl_1) = \int dz\, f^0(z) \, H((x+sv,v)_{s \in [0,T]})$$
and that the equation is complemented with the condition $d\tanaka{t}(\gl) = 0$ for $t < \tau_{n -1}$
when $|\gl| = n$.

\medskip

Formally, the evolution equation \eqref{eq: mild tanaka} is characteristic of a coagulation process
 (see e.g.\,\cite{Norris_clusters}). The limiting particle process behaves as a  coagulation process at the level of the  cluster paths with aggregation rates depending on the law of the process itself.
It would be interesting to derive a martingale interpretation of the coagulation process in the spirit of Tanaka's process \cite{Tanaka} (see also the surveys \cite{Sznitman_saint_flour, Meleardsurvey}), 
which describes the typical evolution of a particle velocity associated with the homogeneous Boltzmann equation.

\bigskip

\noindent 
{\it Proof of Theorem \ref{Prop: convergence dynamical clusters}}

As in the proof of Proposition \ref{Prop: convergence Boltzmann}, 
the concentration estimate \eqref{eq: covariance Pi} implies that 
the convergence in probability \eqref{eq: convergence clusters en proba} 
will simply follow from  the convergence of the expectation $\bbE_\eps \left[  \Pi^\eps_{[0,T]}(H) \right]$
defined in \eqref{eq: esperance tanakae}.
This limit can be deduced by  
Proposition \ref{prop: cluster expansion limit} after a suitable change of variables
\begin{equation}
\label{eq: esperance limite  tanaka}
\lim_{\mu_\eps \to \infty} \bbE_\eps \left[  \Pi^\eps_{[0,T]}(H)   \right]  
= \int   d \tanaka{T} (\gl_1) H(\gl_1) ,
\end{equation}
where $\gl_1$ is indexed by the parameters of the form \eqref{eq: cluster path} and the limiting distribution of a typical cluster path follows from \eqref{eq: def tanakae} 
\begin{align}
\label{eq: def tanaka}
d \tanaka{T} & (\gl_1)  =  d \nu_{T}^{(0)} (\lambda_1) \\
& + \sum_{k \geq 2}\frac{(-1)^{k-1} }{(k-1)!}   \sum_{\cT^{ov} _\prec}  
 d \nu_{T}^{(0)} (\lambda_1) 
\int  d \Theta_{k-1} \; d\Omega_{k-1}  
d\tilde \nu_{T}^{(0)} (\lambda_2) \dots d\tilde \nu_{T}^{(0)} (\lambda_k)
\prod_{ e = \{ i,j\}  \in E(  \cT_\prec ^{ov} )}    \big(  ( {\textbf  v}_i (\tau_e) - {\textbf  v}_j  (\tau_e) ) \cdot\omega_{e }\big)_+ \, \delta_{{\textbf  x}_i (\tau_e ) - {\textbf  x}_j(\tau_e )}  .
\nonumber
\end{align}
The subscript $T$ stresses that the collision times and the overlap times belong to the time interval  $[0,T]$.
Note that the remaining translational degree of freedom is now fixed by the root of $\lambda_1$.
The previous representation remains valid at time $t \in [0,T]$ : the distribution of the process, denoted by $\tanaka{t}$, is obtained by restricting formula \eqref{eq: def tanaka} to the cluster paths formed during the time interval $[0,t]$
\begin{equation}
\label{eq: esperance limite  tanaka}
\lim_{\mu_\eps \to \infty} \bbE_\eps \left[  \Pi^\eps_{[0,t]}(H)   \right]  = \int d \tanaka{t} (\gl_1) H(\gl_1) .
\end{equation}
This completes the proof of the convergence in probability \eqref{eq: convergence clusters en proba}.


\medskip


We are going to derive \eqref{eq: mild tanaka} and show that this limiting structure can be interpreted as a coagulation process at the level of the  cluster paths. We follow a strategy similar to one of 
the proof of Proposition \ref{Prop: convergence Boltzmann} : we first simplify the measure $d \tanaka{T} $ and then identify the coagulation equation.

\medskip

\noindent
\textbf{Step 1} : {\it Restriction to a relevant time ordered path structure.}

\begin{figure}[h] 
\centering
\includegraphics[width=4in]{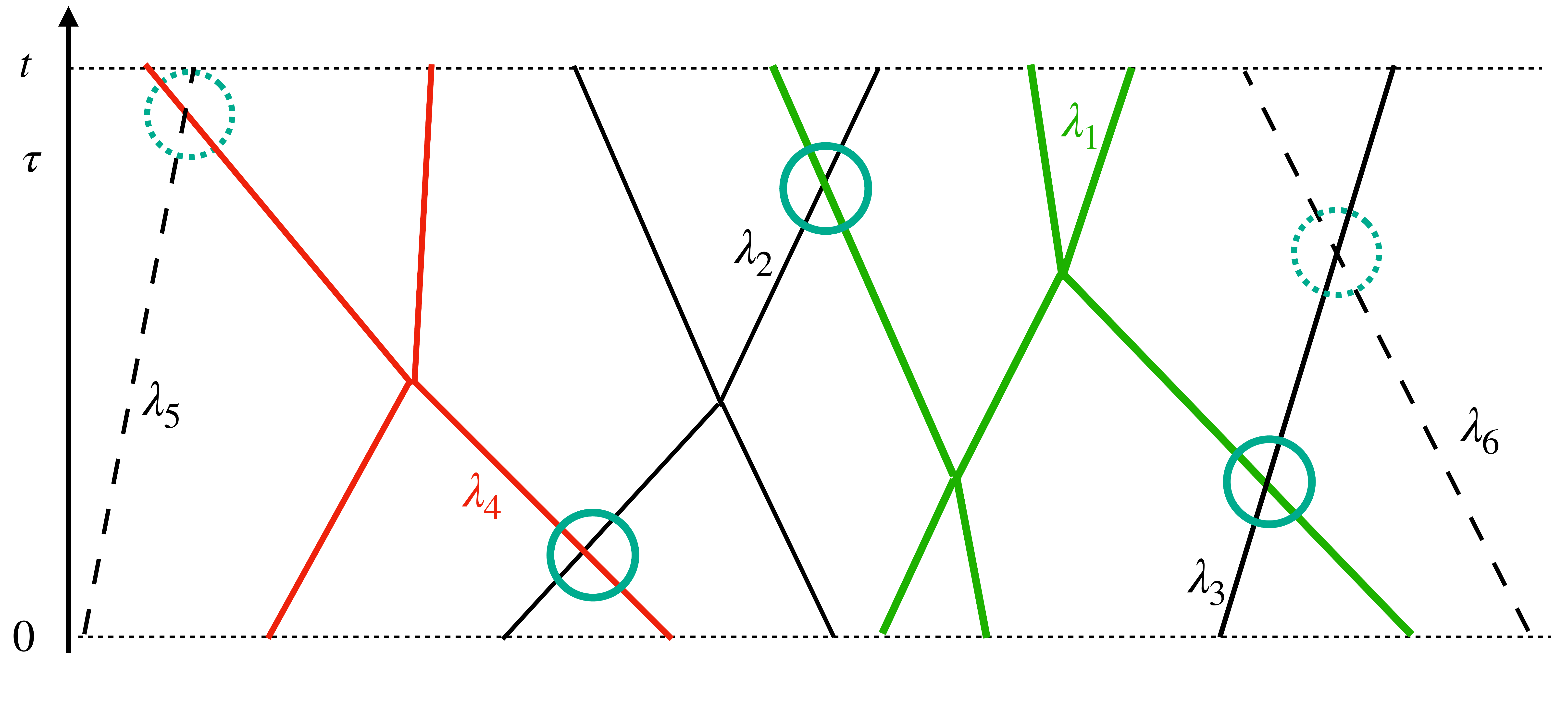} 
\hskip1cm
\includegraphics[width=2in]{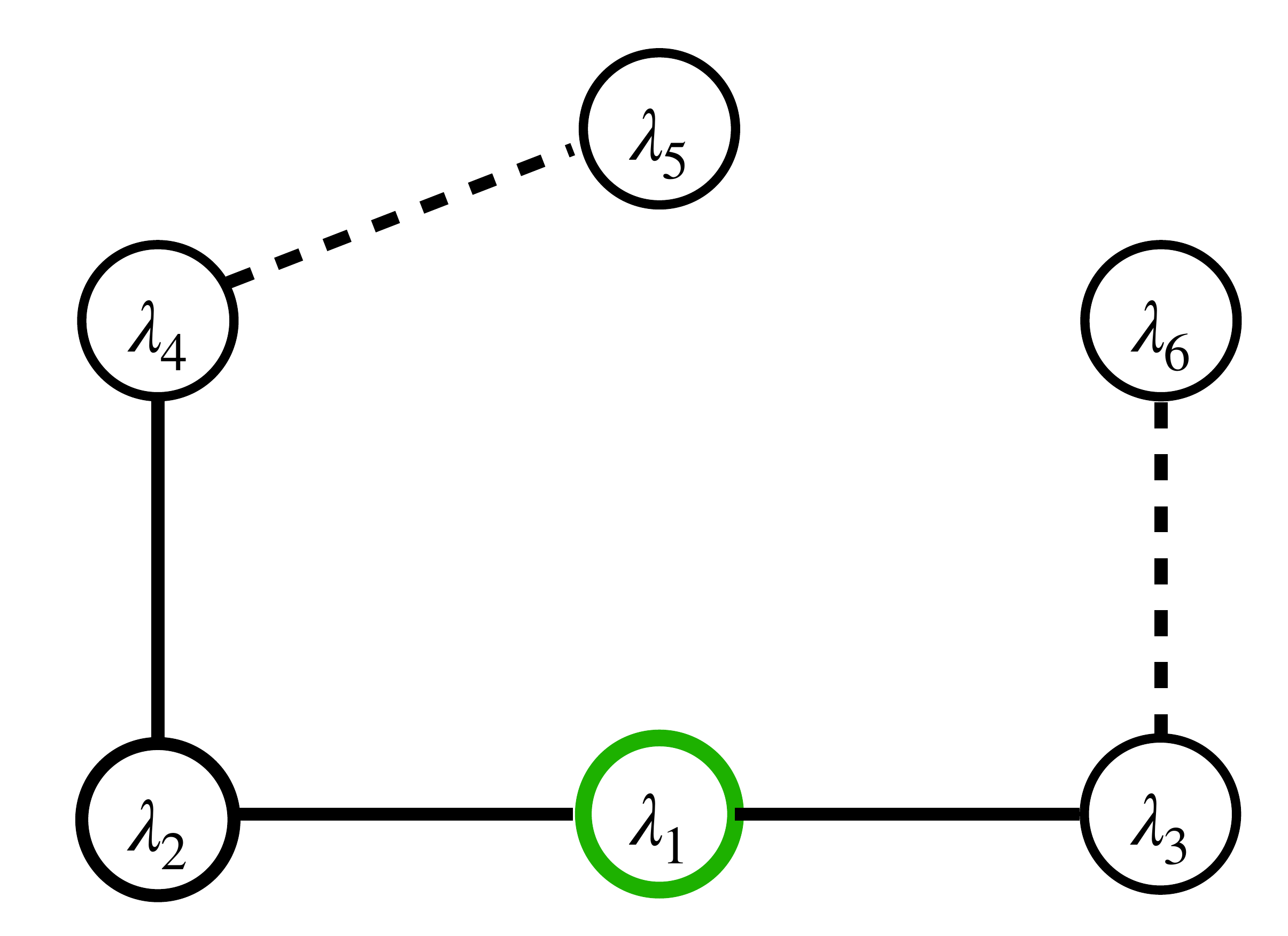} 
\caption{\small On the left, $6$ cluster paths are depicted and
the circles indicate the overlaps between the cluster paths.
The corresponding overlap graph $\cT_\prec^{ov}$ is depicted on the right : the relevant part $\cA$ of the overlap graph $\cT_\prec^{ov}$ is represented by thick lines and the edges which can be neglected by dashed lines.
The edges are associated with the circles representing the overlaps on the left figure.
 The last relevant overlap occurs at time $\tau$ between $\gl_1$ and $\gl_2$.
}
\label{figure: path cluster reduction}
\end{figure}

First, we are going to simplify the series \eqref{eq: def tanaka} by keeping only the relevant cluster overlaps needed to determine the structure of $\gl_1$ (see Figure \ref{figure: path cluster reduction}).
Given $\cT_\prec^{ov}$ with edges ordered according to the overlap times, we build a growing collection of trees $\cA_1 \subset \cA_2 \subset \dots$ as follows.
Starting from the vertex $\cA_1 = \{ 1 \}$ associated with $\gl_1$, all the vertices and edges connected to $1$ are added to form the set $\cA_2$. Suppose that $\ell$ belongs to $\cA_2$ and that the edge 
$\{1,\ell\}$ has order $k$, 
then all the neighbours of $\ell$ are added to $\cA_3$
provided they are linked to $\ell$ by an edge with an order smaller than $k$, i.e.\,if the corresponding overlap has occurred before $k$. Iterating this procedure leads to an ordered tree $\cA$ indexing a set of relevant cluster paths
(see Figure \ref{figure: path cluster reduction}).

\begin{figure}[h] 
\centering
\includegraphics[width=4in]{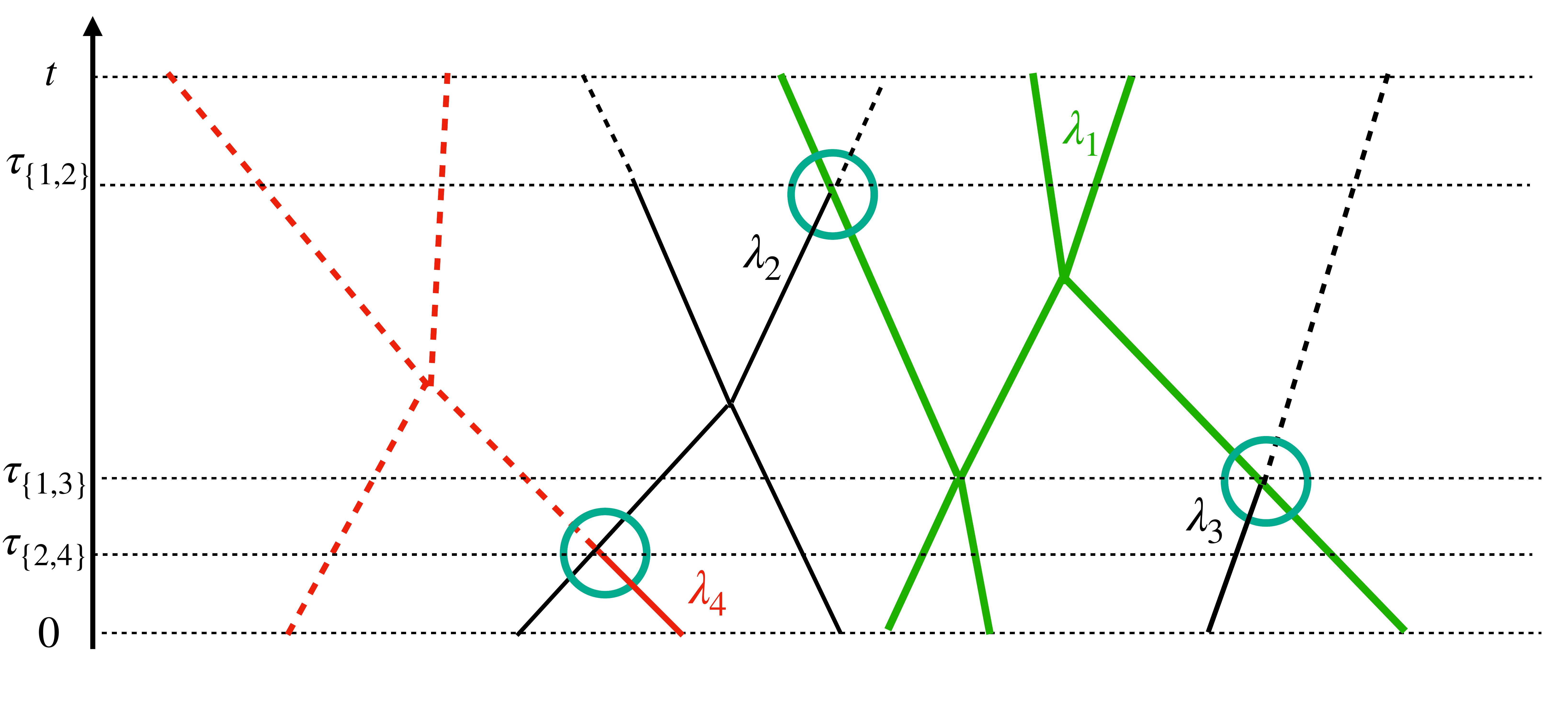} 
\caption{\small After removing the non necessary overlaps as in Figure \ref{figure: path cluster reduction}, the remaining cluster paths have been restricted to the time intervals before the overlap times. 
Note that the number of particles in the cluster path might be reduced by this procedure : in the case of $\gl_4$, depicted above, the branching beyond the overlap time $\tau_{\{2,4\}}$ is discarded so that  the forest $\gl_4$ in the time interval 
$[0,\tau_{\{2,4\}}]$ has only one particle.
Ultimately, the restricted cluster paths, represented with plain lines, are sufficient to compute the  coefficients
in the series decomposition of $\tanaka{t} (\gl_1)$. 
}
\label{figure: nettoyage_total}
\end{figure}

Finally, we consider  the  cluster paths indexed by  $\cA$  only up to the overlap time as explained in Figure \ref{figure: nettoyage_total}. This reduced description is sufficient to compute the terms of the series 
\eqref{eq: def tanaka}. Indeed the branchings and the overlaps which have been discarded  compensate  as in Formula \eqref{eq: compensation}.

By analogy with the proof of Proposition \ref{Prop: convergence Boltzmann}, 
we  denote by   $\pmb{\cA}_k^{\prec}$ the set containing the trees with $k$ vertices of the previous form, i.e.\,the trees rooted in the vertex indexing $\gl_1$ such that the edge orders are decreasing when examined from the root to a leaf.
In summary, to reconstruct $\gl_1$, it is enough to prescribe a tree $\cA$ in $\pmb{\cA}_k^{\prec}$,  an ordered collection $\Theta_{k-1} $ of overlapping times  and $\Omega_{k-1}$, the associated collection of deflection parameters.  The distribution \eqref{eq: def tanaka} of $\gl_1$ can be rewritten for $t \leq T$ as 
\begin{align}
\label{eq: reduction tanaka}
d \tanaka{t} & (\gl_1) = d \nu_{t}^{(0)} (\lambda_1) \\
& + \sum_{k \geq 2} (-1)^{k-1}   \sum_{\cA \in \pmb{\cA}_k^{\prec}}  
d \nu_{t}^{(0)} (\lambda_1) 
\int   d \Theta_{k-1} \; d\Omega_{k-1}  d\tilde \nu_{\tau_2}^{(0)} (\lambda_2) \dots d\tilde \nu_{\tau_k}^{(0)} (\lambda_k)
 \prod_{ e = \{i,j\} \in E(  \cA )}    \big(  ( {\textbf  v}_i (\tau_e) - {\textbf  v}_j  (\tau_e) ) \cdot\omega_{e }\big)_+ \, \delta_{{\textbf  x}_i (\tau_e ) - {\textbf  x}_j(\tau_e )}  .
\nonumber
\end{align}
Note that the factor $1/(k-1)!$ in \eqref{eq: def tanaka} is no longer needed as the cluster paths are indexed by the time ordering.
This is the counterpart of Formula \eqref{eq: limite distribution Boltzmann def} for the particle density. 

\medskip

\noindent
\textbf{Step 2} : {\it Identification of the coagulation equation.}

From \eqref{eq: reduction tanaka}, we are now going to derive \eqref{eq: mild tanaka} which is stated below in an 
 asymmetric way 
\begin{align}
\label{eq: mild tanaka 0}
\int  d \tanaka{t} (\gl_1) H (\gl_1) 
 = & \int d \tanaka{0} (\gl_1) H (\gl_1)   \\
&  + \int  d \tau  \, d \omega \, d \tanaka{\tau} (\gl_1) \, d \tanaka{\tau} (\gl_2) 
\sum_{ i \in \gl_1 \atop j \in \gl_2 }  \Big(    \frac{1}{2} H \big( [\gl_1 \wedge \gl_2]^{i,j,\tau, \omega} \big) -   H (\gl_1) \Big) 
\, \delta_{{\textbf  x}_i (\tau) - {\textbf  x}_j(\tau)}   \,   
\big(  ( {\textbf  v}_i (\tau) - {\textbf  v}_j  (\tau) ) \cdot \omega \big)_+ \;.  \nonumber
\end{align}
Consider the terms in the series \eqref{eq: reduction tanaka} with at least one overlap ($k \geq 2$)
or such that $| \gl _1| \geq 2$. Then starting from time $t$ and looking backward at the  trajectories   of $\gl_1$  and at the overlap times $\Theta_{k-1}$, we denote by $\tau$ the first time at which a collision or an overlap occur.
\begin{itemize}
\item If an overlap occurs at time $\tau$, then the trajectory of $\gl_1$ is unchanged at time $\tau$ so that the restriction of the cluster path  from $[0,t]$ to $[0,\tau]$ is still given by the same set of parameters of the form \eqref{eq: cluster path}.
Splitting the tree $\cA$ (coding the overlaps) into two parts containing respectively $\gl_1$ and $\gl_2$, the cluster paths overlapping $\gl_1$ during $[0,\tau]$ (resp $\gl_2$) can be grouped to reconstruct the 
distributions $\tanaka{\tau} (\gl_1)$ (resp $\tanaka{\tau} (\gl_2)$). 
In this way, the loss term in \eqref{eq: mild tanaka 0} is recovered.
\item If a collision in $\gl_1$ occurs at time $\tau$ between particles $i$ and $j$ with deflection parameter $\omega$, the collision tree $\cT_{\prec}$ associated with $\gl_1$ can be split into two parts $\cT_{\prec,1}$ and $\cT_{\prec,2}$
defining two cluster paths  $\gl_1', \gl_2'$ so that,  by the definition \eqref{eq: merging cluster path},
\begin{equation}
\gl_1 = [\gl_1' \wedge \gl_2']^{i,j,\tau, \omega}.
\end{equation}
Proceeding as before and regrouping the  cluster paths overlapping $\gl_1'$ and $\gl_2'$,  the distributions $\tanaka{\tau} (\gl_1'), \tanaka{\tau} (\gl_2')$ can be identified.
The factor $1/2$ is necessary as the splitting of $\gl_1$ is symmetric.
Changing variables from $\gl_1', \gl_2'$ to $\gl_1, \gl_2$, this leads to the gain term in \eqref{eq: mild tanaka 0}.
\end{itemize}
Proposition  \ref{Prop: convergence dynamical clusters} is therefore complete.

\bibliographystyle{abbrv}
\bibliography{biblio_Boltzmann}

\end{document}